\date{July 11 of 2006}

\documentclass[11pt,british,letterpaper]{amsart}

\usepackage{babel,eucal,url,amssymb,stmaryrd,
enumerate,amscd,
}

\theoremstyle{plain}
\newtheorem{thm}{Theorem}[section]
\newtheorem{lem}[thm]{Lemma}
\newtheorem{pro}[thm]{Proposition}
\newtheorem{co}[thm]{Corollary}

\theoremstyle{definition}
\newtheorem{defn}[thm]{Definition}

\theoremstyle{remark}
\newtheorem{rem}[thm]{Remark}

\numberwithin{equation}{section}

\newcommand{\lie}[1]{\operatorname{\mathfrak{#1}}}

\newcommand{\Gtwo}{\ifmmode{{\rm G}_2}\else{${\rm G}_2$}\fi}
\newcommand{\frg}{\lie{g}}

\newcommand{\LC}{{\nabla^g}}

\newcommand{\Hodge}{\mathord{\mkern1mu *}}

\newcommand{\RR}{{\mathbb{R}}}
\newcommand{\TT}{{\mathbb{T}}}

\def\sideremark#1{\ifvmode\leavevmode\fi\vadjust{\vbox to0pt{\vss
 \hbox to 0pt{\hskip\hsize\hskip1em
 \vbox{\hsize2.5cm\tiny\raggedright\pretolerance10000
 \noindent #1\hfill}\hss}\vbox to8pt{\vfil}\vss}}}%

\begin{document}

\title[Compact Nearly K\"ahler $6$-manifolds with conical singularities]
{Nearly hypo structures and\\
compact Nearly K\"ahler $6$-manifolds with conical singularities}
\date{\today}

\author{Marisa Fern\'andez}
\address[Fern\'andez]{Universidad del Pa\'{\i}s Vasco\\
Facultad de Ciencia y Tecnolog\'{\i}a, Departamento de Matem\'aticas\\
Apartado 644, 48080 Bilbao\\ Spain} \email{marisa.fernandez@ehu.es}

\author{Stefan Ivanov}
\address[Ivanov]{University of Sofia "St. Kl. Ohridski"\\
Faculty of Mathematics and Informatics\\
Blvd. James Bourchier 5\\
1164 Sofia, Bulgaria} \email{ivanovsp@fmi.uni-sofia.bg}

\author{Vicente Mu\~noz}
\address[Mu\~noz]{Instituto de Ciencias Matem\'aticas CSIC-UAM-UC3M-UCM\\
Consejo Superior de Investigaciones Cient\'{\i}ficas\\
Serrano 113bis\\
28006 Madrid, Spain} \email{vicente.munoz@imaff.cfmac.csic.es}

\author{Luis Ugarte}
\address[Ugarte]{
Departamento de Matem\'aticas\,-\,I.U.M.A.\\
Universidad de Zaragoza\\
Campus Plaza San Francisco\\
50009 Zaragoza, Spain} \email{ugarte@unizar.es}

\begin{abstract}
We prove that any totally geodesic hypersurface $N^5$ of a
$6$-dimensional nearly K\"ahler manifold $M^6$ is a Sasaki-Einstein
manifold, and so it has a hypo structure in the sense of
\cite{ConS}. We show that any Sasaki-Einstein 5-manifold defines a
nearly K\"ahler structure on the sin-cone $N^5\times\mathbb R$, and
a compact nearly K\"ahler structure with conical singularities on
$N^5\times [0,\pi]$ when $N^5$ is compact thus providing a link
between Calabi-Yau structure on the cone $N^5\times [0,\pi]$ and the
nearly K\"ahler structure on the sin-cone $N^5\times [0,\pi]$. We
define the notion of {\it nearly hypo} structure that leads to a
general construction of nearly K\"ahler structure on
$N^5\times\mathbb R$. We determine {\it double hypo} structure as
the intersection of hypo and nearly hypo structures and classify
double hypo structures on 5-dimensional Lie algebras with non-zero
first Betti number. An extension of the concept of nearly K\"ahler
structure is introduced, which we refer to as {\it nearly half flat}
$SU(3)$-structure, that leads us to generalize the construction of
nearly parallel $G_2$-structures on $M^6\times\mathbb R$ given in
\cite{BM}. For $N^5=S^5\subset S^6$ and for $N^5=S^2 \times
S^3\subset S^3 \times S^3$, we describe explicitly a Sasaki-Einstein
hypo structure as well as the corresponding nearly K\"ahler
structures on $N^5\times\mathbb R$ and $N^5\times [0,\pi]$, and the
nearly parallel $G_2$-structures on $N^5\times\mathbb R^2$ and
$(N^5\times [0,\pi])\times [0,\pi]$.
\end{abstract}

\maketitle

\setcounter{tocdepth}{2} \tableofcontents

\section{Introduction}

Let $N^5$ be a $5$-manifold with an $SU(2)$-structure, that is, the
frame bundle of $N^5$ has a reduction to the group $SU(2)$.
Recently, Conti and Salamon \cite{ConS} have proved that such a
structure is determined  by a quadruplet
$(\eta,\omega_1,\omega_2,\omega_3)$ of differential forms, that we
shall abbreviate as $(\eta,\omega_i)$, where $\eta$ is a $1$-form
and $\omega_i$ are $2$-forms satisfying certain relations (see
Section \ref{hypo}). An $SU(2)$-structure $(\eta,\omega_i)$ is said
to be {\em hypo} if the $2$-form $\omega_1$ and the $3$-forms
$\eta\wedge\omega_2$ and $\eta\wedge\omega_3$ are closed.

Hypo geometry is a generalization of Sasaki-Einstein geometry. In
fact, any Sasaki-Einstein $5$-manifold has an $SU(2)$-structure
$(\eta,\omega_i)$, where $\eta$ is the contact form, that satisfies
the differential equations
  \begin{equation}\label{e-s}
  d\eta=-2\omega_3,  \quad d\omega_1=3\eta\wedge\omega_2,   \quad
  d\omega_2=-3\eta\wedge\omega_1,
  \end{equation}
and so is a hypo structure, after interchanging the form $\omega_1$
with $\omega_3$. This is due to the following. A Sasaki-Einstein
$5$-manifold $N^5$ may be defined as a Riemannian manifold such that
$N^5\times\mathbb R$ with the cone metric is K\"ahler and Ricci flat
\cite{BGal}, that is, it has  holonomy contained in $SU(3)$ or,
equivalently, its $SU(3)$-structure is integrable. This means that
there is an almost Hermitian structure, with K\"ahler form $F$, and
a complex volume form $\Psi = \Psi_{+} + i \Psi_{-}$ on
$N^5\times\mathbb R$ satisfying $dF=d\Psi_+=d\Psi_-=0$. But an
integrable $SU(3)$-structure on the cone $N^5\times\mathbb R$
induces an $SU(2)$-structure on $N^5$ satisfying \eqref{e-s} (see
Section~\ref{hypo} for details).

Our goal in this paper is twofold: on the one hand, to show that
Sasaki-Einstein (hypo) $5$-manifolds are  also closely related with
nearly K\"ahler $6$-manifolds (weak holonomy
$SU(3)$ manifolds) giving a method to construct nearly K\"ahler
manifolds from Sasaki-Einstein $5$-manifolds; and on the other hand,
to give a method of construction of nearly parallel $G_2$-structures
on $M^6\times\mathbb R$ starting from certain $SU(3)$-structures on
$M^6$,  which we call {\em nearly half flat}, leading to a
generalization of the construction given in~\cite{BM}.

To this end, in Section \ref{hypo} it is shown that any totally
geodesic hypersurface $N^5$ of a nearly K\"ahler $6$-manifold $M^6$
has a natural Sasaki-Einstein $SU(2)$-structure
$(\eta,\omega_1,\omega_2,\omega_3)$ satisfying~\eqref{e-s}.
Furthermore, the converse also holds. In fact, we prove that any
Sasaki-Einstein $SU(2)$-structure on $N^5$ satisfying (\ref{e-s})
defines an $SU(3)$-structure on the sin-cone $N^5\times\mathbb R$
which is nearly K\"ahler (see Theorem \ref{newnk} in Section
\ref{nearlyparallel}). Actually, our result is slightly more general
and it applies to {\it nearly hypo} $SU(2)$-structures satisfying
the {\em evolution nearly hypo equations} established  in
Proposition \ref{evpro}. Nearly hypo structures are the natural
$SU(2)$-structures induced on oriented hypersurfaces of nearly
K\"ahler $6$-manifolds. In particular, when $N^5$ is a compact
Sasaki-Einstein $SU(2)$-manifold, one gets a compact nearly K\"ahler
structure with conical singularities on $N^5\times [0,\pi]$.

Returning to a Sasaki-Einstein structure, it can be defined as a
structure whose cone is K\"ahler and Ricci flat. We show (see
Corollary \ref{conewnk}) that in dimension $5$ a Sasaki-Einstein
structure could also be defined as a structure whose sin-cone is
nearly K\"ahler (weak holonomy $SU(3)$ manifold). In this way,
Sasaki-Einstein $5$-manifolds provide a link between Calabi-Yau
cones and nearly K\"ahler (weak holonomy $SU(3)$) sin-cones.

More general, in Section \ref{double-hypo-structures}, we define a
\emph{double hypo structure} as an $SU(2)$-structure which is hypo
and nearly hypo; a diagram representing the relations among the
classes of $SU(2)$-structures is inserted. In Section
\ref{double-hypo-structures}, we show also that double hypo
structures are precisely those $SU(2)$-structures whose sin-cone
carry a half-flat $SU(3)$-structure. We describe all 5-dimensional
Lie algebras with non-zero first Betti number which have a double
hypo structure, and prove that solvable Lie groups cannot admit
invariant double hypo structures. Double hypo structures give a
relation between Calabi-Yau solution to Conti-Salamon evolution
equations \eqref{vol} and nearly K\"ahler solution to the nearly
hypo evolution equations \eqref{evolunk} discovered in
Proposition~\ref{evpro} below.

In \cite{BM} it is proved that if $M^6$ is a nearly K\"ahler
manifold, then the sin-cone $M^6\times\mathbb R$ has a natural
nearly parallel $G_2$-structure. We generalize  this construction of
nearly parallel $G_2$-structures proving (see Proposition
\ref{g2evpro}) that any {\em nearly half flat} $SU(3)$-structure
$(F,\Psi_{+},\Psi_{-})$ on $M^6$, which means that
$d\Psi_{-}=-2F\wedge F$, can be lifted to a nearly parallel
$G_2$-structure on $M^6\times\mathbb R$ if and only if it satisfies
the {\em evolution nearly half flat equation} (\ref{evolug2})
established in Section~\ref{nearlyhalf}. In this section we insert
two figures, one of them represents the relations among
the classes of $SU(3)$-structures, and the other illustrates how it is possible to
get special $G_2$-metrics by evolution from $SU(3)$-structures.

In Section \ref{examples} we  consider the oriented hypersurfaces
$N^5=S^5\subset S^6$ and $N^5=S^2 \times S^3\subset S^3 \times S^3$.
Since $S^5\subset S^6$ is totally geodesic in $S^6$ with the metric
of the nearly K\"ahler structure on $S^6$, it induces a
Sasaki-Einstein hypo structure on  $S^5$ satisfying \eqref{e-s}. We
describe explicitly such a structure on $S^5$ as well as the nearly
K\"ahler structure on $S^5\times\mathbb R$  and the nearly parallel
$G_2$-structure on $S^5\times\mathbb R^2$.

For $S^2 \times S^3\subset S^3 \times S^3$ we notice that $S^2
\times S^3$ is not totally geodesic in $S^3 \times S^3$ with the
metric of the nearly K\"ahler structure, and we see that the
$SU(2)$-structure induced on $S^2 \times S^3$ is hypo but it does
not satisfy the first equation of \eqref{e-s}. We modify it a little
to obtain a Sasaki-Einstein $SU(2)$-structure on $S^2 \times S^3$
satisfying equations \eqref{e-s}, and then we describe the nearly
K\"ahler structure on $S^2 \times S^3\times\mathbb R$  and the
nearly parallel $G_2$-structure on $S^2 \times S^3\times\mathbb
R^2$.

Finally, we use the recently discovered in \cite{GMS} infinite
family of explicit compact Sasaki-Einstein 5-manifold $Y^{p,q}$ to
construct infinite family of compact nearly K\"ahler manifold with
conical singularities on $Y^{p,q}\times[0,\pi]$.

\section{Hypo structures on $5$-manifolds}\label{hypo}

In this section we show  that any totally geodesic  hypersurface of
a nearly K\"ahler manifold has a Sasaki-Einstein $SU(2)$-structure
satisfying \eqref{e-s}. First we need to recall some properties of
$SU(2)$-structures and, in particular, of hypo structures on
$5$-manifolds.

Consider a $5$-manifold $N^5$ with an $SU(2)$-structure
$(\eta,\omega_1,\omega_2,\omega_3)$, that is to say, $\eta$ is a
$1$-form and $\omega_i$ are $2$-forms on $M$ satisfying
  \begin{equation}\label{defsu2}
  \omega_i\wedge\omega_j=\delta_{ij}v, \quad
  v\wedge\eta\not=0,
  \end{equation}
for some $4$-form $v$, and
  \begin{equation}\label{defsu2-2}
  X\lrcorner\omega_1=Y\lrcorner\omega_2\Rightarrow \omega_3(X,Y)\ge 0,
  \end{equation}
where $X\lrcorner$ denotes the contraction by $X$. Then, it induces
an $SU(3)$-structure $(F,\Psi_+,\Psi_-)$ on $N^5\times\mathbb R$
defined by
  \begin{equation}\label{hyp0}
  F=\omega_1+\eta\wedge dt, \quad
  \Psi=\Psi_++i\Psi_-=(\omega_2+i\omega_3)\wedge(\eta+idt),
  \end{equation}
where $t$ is a coordinate on $\mathbb R$.

Vice versa, let $f: N^5\longrightarrow M^6$ be an oriented
hypersurface of a $6$-manifold $M^6$ with an $SU(3)$-structure
$(F,\Psi_+,\Psi_-)$, and denote by $\mathbb N$ the unit normal
vector field. Then the $SU(3)$-structure induces an
$SU(2)$-structure $(\eta,\omega_1,\omega_2,\omega_3)$ on $N^5$
defined by the equalities~\cite{ConS}
  \begin{equation}\label{hyp1}
  \eta=-\mathbb N\lrcorner F,\quad
  \omega_1=f^*F,\quad \omega_2=\mathbb N\lrcorner\Psi_-, \quad
  \omega_3=-\mathbb N\lrcorner\Psi_+.
  \end{equation}
An $SU(2)$-structure determined by $(\eta,\omega_i)$ is called {\em
hypo} if it satisfies the equations \cite{ConS}
  \begin{equation}\label{rhypo}
  d\omega_1=0, \qquad d(\eta\wedge\omega_2)=0,
  \qquad d(\eta\wedge\omega_3)=0.
  \end{equation}

Suppose that $M^6$ has holonomy contained in $SU(3)$, that is, the
$SU(3)$-structure $(F,\Psi_+,\Psi_-)$ is integrable (i.e. Calabi-Yau) or,
equivalently,
  $$
  dF=d\Psi_+=d\Psi_-=0.
  $$
It is not hard to see that any oriented hypersurface $N^5$ of $M^6$
is naturally endowed with a hypo structure \cite{ConS}. Indeed, the
conditions $dF=d\Psi_+=d\Psi_-=0$ imply that the induced
$SU(2)$-structure on $N^5$ defined by \eqref{hyp1} satisfies
\eqref{rhypo}. Regarding the converse, Conti and Salamon \cite{ConS}
prove that a real analytic hypo structure on $N^5$ (that is, when
$N^5$ and the reduction of the frame bundle of $N^5$ both are
analytic) can be lifted to an integrable $SU(3)$-structure on
$N^5\times \mathbb R$, that is, $(\eta,\omega_i)$ belongs to a
one-parameter family of hypo structures $(\eta(t),\omega_i(t))$
satisfying the evolution equations
  \begin{equation}\label{vol}
  \left\{\begin{array}{l}
  \partial_t\omega_1=-d\eta\\
  \partial_t(\eta\wedge\omega_3)=d\omega_2\\
   \partial_t(\eta\wedge\omega_2)=-d\omega_3.
  \end{array} \right.
  \end{equation}

Next we study totally geodesic hypersurfaces of nearly K\"ahler
$6$-manifolds $M^6$, that is, $M^6$ has an $SU(3)$-structure
$(F,\Psi_+,\Psi_-)$ which satisfies the following differential
equations \cite{Hit}
  \begin{equation}\label{nkdef}
  dF=3\Psi_+, \qquad d\Psi_-=-2F\wedge F.
  \end{equation}

\begin{lem}\label{lem1}
If $f: N^5\longrightarrow M^6$ is a totally geodesic hypersurface of
a nearly K\"ahler manifold $M^6$, then the induced $SU(2)$-structure
\eqref{hyp1} on $N^5$ satisfies the differential equations
\eqref{e-s}.
\end{lem}

\begin{proof}
Let $(M^6,g,F,\Psi_+,\Psi_-)$ be a nearly K\"ahler $6$-manifold. The
Nijenhuis tensor $N$ is a $3$-form $N=-\Psi_-$ and it is parallel
with respect to the Gray characteristic connection $\nabla$
\cite{Kir}. This connection was defined by Gray \cite{Gr2,Gr3,Gr1}
and it turns out to be the  unique linear connection preserving the
nearly K\"ahler structure and having totally skew-symmetric torsion
$T=N=-\Psi_-$ \cite{FI1}, i.e.
  \begin{equation}\label{nkcon}
  \nabla=\LC +\frac12 T=\LC-\frac12\Psi_-, \qquad
  \nabla\Psi_-=0,
  \end{equation}
where $\LC$ is the Levi-Civita connection of the metric $g$.

We calculate using \eqref{hyp1} and \eqref{nkdef} that
  \begin{gather}\label{sas2}
  d\omega_1=d(f^*F)=3f^*\Psi_+=3\eta\wedge\omega_2, \\\label{sas1}
  d\eta=-d(\mathbb N\lrcorner F)=-(L_{\mathbb N}F)+\mathbb N\lrcorner
  dF = -(L_{\mathbb N}F)-3\omega_3,
  \end{gather}
where $L$ denotes the Lie derivative.

Further, $-(L_{\mathbb N}F)=-\LC_{\mathbb N}F$, since $N^5$ is
totally geodesic. Apply \eqref{nkcon} to the latter equality, take
into account $\nabla F=0$ and \eqref{hyp1} to derive
  \begin{equation}\label{sas11}
  -\LC_{\mathbb N}F= -\nabla_{\mathbb
  N}F-\frac12\sum_{i=1}^6e_i\lrcorner F\wedge e_i\lrcorner(\mathbb
  N\lrcorner\Psi_-)=-\mathbb N\lrcorner\Psi_+=\omega_3,
  \end{equation}
where $\{e_1,\dots, e_6=\mathbb N\}$ is an $SU(3)$ adapted basis.
Substitute \eqref{sas11} into \eqref{sas1} to get the first equality
in \eqref{e-s}.

In view of \eqref{sas2}, it remains to prove the third equality in
\eqref{e-s}. Similarly as above, applying \eqref{hyp1},
\eqref{nkcon} and \eqref{nkdef}, we calculate
  \begin{equation}\label{sas3}
  \begin{aligned}
  d\omega_2&=d(\mathbb N\lrcorner \Psi_-)=
  L_{\mathbb N}\Psi_--\mathbb N\lrcorner d\Psi_-= \LC_{\mathbb
  N}\Psi_-+2\mathbb N\lrcorner(F\wedge F)=\\ &= -\frac14\mathbb
  N\lrcorner(\sum_{j=1}^6e_j\lrcorner\Psi_-\wedge
  e_j\lrcorner\Psi_-)+2\mathbb N\lrcorner(F\wedge F)=\\ &=
  -\frac12\mathbb N\lrcorner(F\wedge F)+2\mathbb N\lrcorner(F\wedge
  F)=\frac32\mathbb N\lrcorner(F\wedge F)=-3\eta\wedge\omega_1,
  \end{aligned}
  \end{equation}
where we have used the identity
  $$
  \sum_{j=1}^6e_j\lrcorner T\wedge e_j\lrcorner T=2F\wedge F
  $$
valid on any nearly K\"ahler $6$-manifold \cite{FI1}.
\end{proof}

\begin{thm}
Any totally geodesic hypersurface $N^5$ of a nearly K\"ahler
$6$-manifold $M^6$ admits a Sasaki-Einstein  hypo structure, and
therefore the Conti-Salamon evolution equations~\eqref{vol} can be
solved for $N^5\times \mathbb R$.
\end{thm}

\begin{proof}
Clearly Lemma~\ref{lem1} implies that the induced $SU(2)$-structure
satisfies \eqref{rhypo}, i.e, it is a hypo structure. Moreover,
Lemma~\ref{lem1} shows that the induced almost contact metric
structure $(\eta,\omega_3)$ on $N^5$ is Sasaki-Einstein. Indeed,
\eqref{e-s} implies that the conical $SU(3)$-structure on
$M=N^5\times\mathbb R$ defined by
  \begin{gather}\label{integrable}
  F=t^2\omega_3+t\eta\wedge dt,   \quad\
  \Psi=t^2(\omega_2+i\omega_1)\wedge(t\eta+idt)
  \end{gather}
satisfies $dF=d\Psi=0$, i.e. it is an integrable $SU(3)$-structure
(see e.g. \cite{BGal}) which clearly is a solution to the
Conti-Salamon evolution equations \eqref{vol}.
\end{proof}

\begin{rem}
We notice  that any Sasaki-Einstein $5$-manifold has a hypo
$SU(2)$-structure which satisfies \eqref{e-s}. In fact, we know that
a Sasaki-Einstein $5$-manifold $N^5$ is such that the cone
$N^5\times \mathbb R$ is K\"ahler and Ricci flat, that is, its
$SU(3)$-structure is integrable, and so induces an $SU(2)$-structure
on $N^5$ satisfying \eqref{e-s} which is equivalent to equations
$(14)$ in \cite{ConS}, although the two forms $\omega_2,\omega_3$
are not given explicitly there since the $SU(3)$-structure on the
cone is not explicit;  we just know that such a structure does exist
and is given by \eqref{integrable}.
\end{rem}

\section{Nearly hypo structures}\label{nearlyparallel}

Let $(\eta, \omega_i)$ be an $SU(2)$-structure on $N^5$ and consider
the $SU(3)$-structure $(F,\Psi_+,\Psi_-)$ on $N^5\times \mathbb R$
defined by \eqref{hyp0}.

We look for sufficient conditions imposed  on the $SU(2)$-structure
$(\eta, \omega_i)$ which imply that the induced $SU(3)$-structure on
$N^5\times \mathbb R$ is nearly K\"ahler, i.e. it satisfies
\eqref{nkdef}.

\begin{defn}
We call an $SU(2)$-structure $(\eta, \omega_i)$ on a $5$-manifold
$N^5$ a \emph{nearly hypo structure} if it satisfies the following
two equations:
  \begin{equation}\label{nkhypo}
  d\omega_1=3\eta\wedge\omega_2,\quad
  d(\eta\wedge\omega_3)=-2\omega_1\wedge\omega_1.
  \end{equation}
\end{defn}

Consider $SU(2)$-structures $(\eta(t),\omega_i(t))$ on $N^5$
depending on a real parameter $t\in\mathbb R$, and the corresponding
$SU(3)$-structures $(F(t),\Psi_+(t),\Psi_-(t))$ on $N^5\times\mathbb
R$. We have

\begin{pro}\label{evpro}
An $SU(2)$-structure $(\eta,\omega_i)$ on $N^5$ can be lifted to a
nearly K\"ahler structure $(F(t),\Psi_+(t),\Psi_-(t))$ on
$N^5\times\mathbb R$ defined by \eqref{hyp0} if and only if it is a
nearly hypo structure which generates an 1-parameter family of
$SU(2)$-structures $(\eta(t),\omega_i(t))$ satisfying the following
\emph{evolution nearly hypo equations}
  \begin{equation}\label{evolunk}
  \left\{ \begin{array}{l}
  \partial_t\omega_1=-d\eta-3\omega_3,\\
  \partial_t(\eta\wedge\omega_3)=d\omega_2+ 4\eta\wedge\omega_1,\\
  \partial_t(\eta\wedge\omega_2)=-d\omega_3.
  \end{array} \right.
  \end{equation}
\end{pro}
\begin{proof}
Take the exterior derivatives in \eqref{hyp0} to get that the
equations \eqref{nkdef} hold precisely when \eqref{nkhypo} and the
first two equalities in \eqref{evolunk} are fulfilled.

It remains to show that the equations \eqref{evolunk} imply that
\eqref{nkhypo} hold for each $t$. Indeed, using \eqref{evolunk}, we
calculate
\begin{gather*}\partial_t(d\omega_1-3\eta\wedge\omega_2)=-3(d\omega_3+\partial_t(\eta\wedge\omega_2) )=0.
\end{gather*}
Hence, the first equality in \eqref{nkhypo} is independent on $t$
and therefore is valid for all $t$ since it holds in the beginning
for $t=0$. Further, using the already proved first equality in
\eqref{nkhypo} as well as the defining equalities \eqref{defsu2}, we
obtain
\begin{gather*}
\partial_t[d(\eta\wedge\omega_3)+2\omega_1\wedge\omega_1]=-4\eta\wedge d\omega_1=0.
\end{gather*}
Hence, both equalities in \eqref{nkhypo} survive in time.
\end{proof}
\begin{rem}
The assumption $(\eta(t),\omega_i(t))$ to be an $SU(2)$-structure
for all t in Proposition~\ref{evpro} can not be avoided as it is
shown in the example described in the last Section~\ref{last}.
\end{rem}

\begin{pro}
Any  $SU(2)$-structure satisfying the two first equations of
\eqref{e-s} is a nearly hypo structure.
\end{pro}
\begin{proof}
The two first equations of \eqref{e-s} together with \eqref{defsu2}
yield
  $$
  d\omega_1=3\eta\wedge\omega_2,\quad
  d(\eta\wedge\omega_3)=-2\omega_3\wedge\omega_3=-2\omega_1\wedge\omega_1.
  $$
\end{proof}
More generally, we have
\begin{pro}\label{hypnew}
Let $f: N^5\longrightarrow M^6$ be an immersion of an oriented
$5$-manifold into a $6$-manifold with a nearly K\"ahler structure.
Then the $SU(2)$-structure induced on $N^5$ is a nearly hypo
structure.
\end{pro}

\begin{proof}
It follows from \eqref{hyp1} that \cite{ConS}
  $$
  \eta\wedge\omega_2=f^*\Psi_+, \quad
  \eta\wedge\omega_3 =f^*\Psi_-.
  $$
Since $f^*$ commutes with $d$, the above equality together with
\eqref{hyp1} and \eqref{nkdef} imply \eqref{nkhypo}.
\end{proof}

Now, a question remains.

{\bf Question 1.} Does the converse of 
Proposition~\ref{hypnew}
hold?, i.e. is it true that any (real analytic) nearly hypo
structure on $N^5$ can be lifted to a nearly K\"ahler structure on
$N^5\times\RR$?

\begin{rem}\label{answer-Q1}
The affirmative answer to this question is equivalent to showing the
existence of a solution of the evolution nearly hypo equations
\eqref{evolunk}. From a private communication with D. Conti
\cite{Conti}, we know that the answer of Question~1 is affirmative,
at least locally, for real analytic nearly hypo structures. In fact,
if $N^5$ is compact, there is a solution to the nearly hypo
evolution equations on $N^5$, i.e. the real analytic nearly hypo
structure on $N^5$ can be lifted to a nearly K\"ahler structure on
$N^5\times I$, for a sufficiently small interval $I$; and if $N^5$
is non-compact, one always has a local solution to these equations,
that is, there is an open set $U\subset N^5$ such that the real
analytic nearly hypo structure on the $5$-manifold $N^5$ can be
lifted to a nearly K\"ahler structure on $U\times I$, for a
sufficiently small interval $I$.
\end{rem}

Now, we prove the main result in this section solving explicitly the
equations \eqref{evolunk} for Sasaki-Einstein $5$-manifolds.

\begin{thm}\label{newnk}
Let $(N^5,\eta,\omega_i)$ be a Sasaki-Einstein $SU(2)$-structure
satisfying \eqref{e-s}. Then the $SU(3)$-structure $F,\Psi_+,\Psi_-$
on $N^5\times\mathbb R$ defined for $0\le t\le \pi$ by
  \begin{equation}\label{newnks}
  \begin{aligned}
  & F =\sin^2t\left(\sin t\, \omega_1+\cos t\, \omega_3\right)+\sin t\,
  \eta\wedge dt,\\ & \Psi_+ =\sin^3t\, \eta\wedge\omega_2 -
  \sin^2t\left(-\cos t\, \omega_1+\sin t\, \omega_3\right)\wedge dt,
  \\
  & \Psi_- =\sin^3t\left(-\cos t\, \omega_1+\sin t\,
  \omega_3\right)\wedge\eta + \sin^2t\, \omega_2\wedge dt,
  \end{aligned}
  \end{equation}
is a nearly K\"ahler structure on $N^5\times\mathbb R$ generating
the well known Einstein metric
  $$
  g_6=dt^2+\sin^2t\, g_5,
  $$
where $g_5$ is the Sasaki-Einstein metric on $N^5$.

If $(N^5,\eta,\omega_i)$ is compact then $(N^5\times
[0,\pi],F,\Psi_+,\Psi_-)$ is a compact nearly K\"ahler $6$-manifold
with two conical singularities at $t=0$ and  $t=\pi$.
\end{thm}

\begin{proof}
Consider the  $SU(2)$-structure $(\eta(t),\omega_i(t))$ depending on
a real parameter $t$:
  \begin{equation}\label{evonk}
  \begin{aligned}
  &\eta(t)=\sin t\, \eta, \\
  &\omega_1(t)=\sin^2t\left(\sin t\,\omega_1+\cos
  t\,\omega_3\right),\\ & \omega_2(t)=\sin^2t\,\omega_2, \\
  &\omega_3(t)=\sin^2t\left(-\cos t\,\omega_1+\sin
  t\,\omega_3\right).
  \end{aligned}
  \end{equation}
Applying \eqref{e-s} and \eqref{defsu2}, we see that the structure
defined by \eqref{evonk} satisfies the nearly hypo structure
conditions \eqref{nkhypo} as well as the nearly hypo evolution
equations \eqref{evolunk}. Consequently, \eqref{newnks} satisfies
\eqref{nkdef} and therefore it is a nearly K\"ahler structure
on~$N^5\times\mathbb R$.
\end{proof}

As a consequence of the proof of Theorem \ref{newnk}, we derive
\begin{co}\label{conewnk}
An $SU(2)$-manifold $(N^5,\eta,\omega_i)$ is Sasaki-Einstein if and
only if the sin-cone $(N^5\times\mathbb R,F,\Psi_+,\Psi_-)$ with the
$SU(3)$-structure defined by \eqref{newnks} is a nearly K\"ahler
manifold for any $0< t<\pi$.
\end{co}
\begin{proof}
The equations \eqref{newnks} imply
$$dF=\sin tdt\wedge[3\sin t\cos t\omega_1- 3\sin^2 t\omega_3 +2\omega_3 +d\eta]+
\sin^2 t(\sin td\omega_1+\cos td\omega_3).
$$
Consequently, $dF=3\Psi_+ \Leftrightarrow
d\omega_1=3\eta\wedge\omega_2, \quad d\eta=-2\omega_3.$ Using this
equivalence, we obtain
$$d\Psi_-+2F\wedge F=\sin^3t[\sin t\omega_3\wedge(d\eta+2\omega_3)-\cos t\omega_1\wedge d\eta]
+\sin^2t(3\omega_1\wedge\eta+d\omega_2)\wedge dt.
$$
Hence, $d\Psi_-=-2F\wedge F \Leftrightarrow
d\omega_2=-3\eta\wedge\omega_1$. Thus, \eqref{e-s} are equivalent to
\eqref{nkdef}  and the proof is complete.
\end{proof}
 \begin{rem}
 Any Sasaki-Einstein $5$-manifold generates, on one hand, a Calabi-Yau
structure on the cone and, on the other hand, it generates a nearly
K\"ahler (weak holonomy $SU(3)$) structure on the sin-cone, thus
giving a link between these two structures in dimension six.
Moreover, Lemma \ref{lem1} shows that any totally geodesic
hypersurface of a nearly K\"ahler $6$-manifold carries a natural
Sasaki-Einstein structure and therefore one gets a non-compact
Calabi-Yau cone genetared by that structure. Vice versa, any totally
umbilic hypersurface in a Calabi-Yau $6$-manifold with shape
operator $A=id$ carries a natural Sasaki-Einstein structure which
could be lifted to a nearly K\"ahler structure on the sin-cone
according to Theorem \ref{newnk}. It seems that a (local)
description of totally geodesic hypersurfaces of a nearly K\"ahler
$6$-manifold as well as the (local) description of totally umbilic
hypersurfaces of a Calabi-Yau $6$-manifold with shape operator equal
to the identity will provide an explicit relation between Calabi-Yau
$6$-manifolds and nearly K\"ahler $6$-manifolds.
\end{rem}

\begin{rem}
There exist nonhomogeneous examples of Sasaki-Einstein $5$-manifolds;
for instance, there are known $14$ nonhomogeneous Sasaki-Einstein
metrics on $S^2\times S^3$ \cite{BG,BGN}. Using these structures we
obtain examples of local nonhomogeneous nearly K\"ahler
$6$-manifolds constructed according to Theorem \ref{newnk} .
\end{rem}

\section{Double hypo structures}\label{double-hypo-structures}
In this section we are interested in the class of $SU(2)$-structures
on a $5$-manifold which are in the intersection class of hypo and
nearly hypo structures.

\begin{defn}
An  $SU(2)$-structure $(\eta,\omega_i)$ on a $5$-manifold is said to
be \emph{double hypo} if it is hypo and nearly hypo simultaneously.
\end{defn}

The following picture illustrates the various classes of
$SU(2)$-structures.

\setlength{\unitlength}{2.5pt}
\begin{picture}(120,70)(-82,-18)
\put(-15,0){\linethickness{1pt}\line(1,0){31}}
\put(-15,0){\linethickness{1pt}\line(0,1){10}}
\put(-15,10){\linethickness{1pt}\line(1,0){31}}
\put(16,0){\linethickness{1pt}\line(0,1){10}}
\put(-14,3){Sasaki-Einstein}
\put(-25,-15){\linethickness{1pt}\line(1,0){90}}
\put(-25,-15){\linethickness{1pt}\line(0,1){35}}
\put(-25,20){\linethickness{1pt}\line(1,0){90}}
\put(65,-15){\linethickness{1pt}\line(0,1){35}}
\put(-11,13.5){double hypo}
\put(-65,-4){\linethickness{1pt}\line(1,0){90}}
\put(-65,-4){\linethickness{1pt}\line(0,1){35}}
\put(-65,31){\linethickness{1pt}\line(1,0){90}}
\put(25,-4){\linethickness{1pt}\line(0,1){35}}
\put(-51,13){hypo}
\put(33,1){nearly hypo}
\put(-75,-20){\linethickness{1pt}\line(1,0){150}}
\put(-75,-20){\linethickness{1pt}\line(0,1){68}}
\put(-75,48){\linethickness{1pt}\line(1,0){150}}
\put(75,-20){\linethickness{1pt}\line(0,1){68}}
\put(-30,38){5-manifolds with $SU(2)$-structure}
\end{picture}

\bigskip

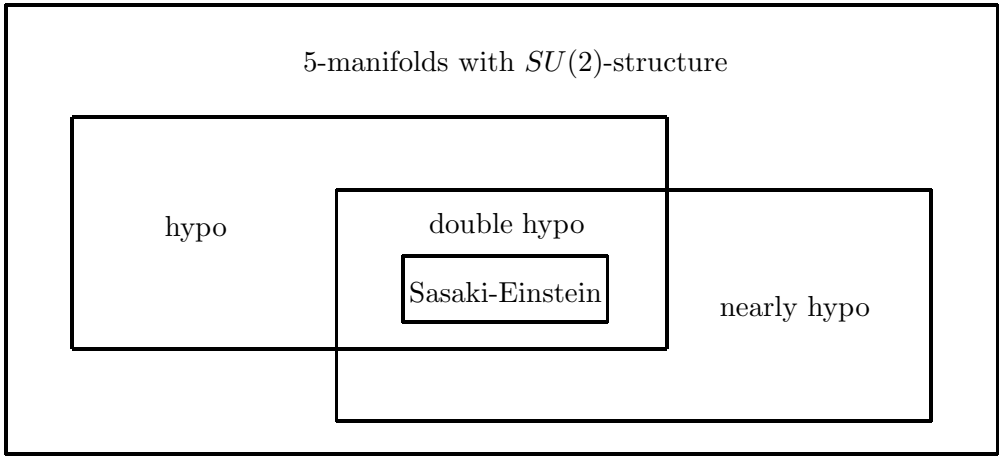
\centerline{{\bf Figure 1}: Classes of $SU(2)$-structures}

\bigskip

Double hypo structures can be lifted in the analytic case, on one
hand, to an integrable $SU(3)$-structure due to the Conti-Salamon
result \cite{ConS} and, on the other hand, taking account of
Remark~\ref{answer-Q1}, to a nearly K\"ahler structure, which
provides a relationship between these distinguished classes of
6-dimensional manifolds:

\setlength{\unitlength}{2.5pt}
\begin{picture}(120,70)(-82,-18)
\multiput(-2,2)(2,0){3}{\circle*{.2}}
\multiput(-7,1)(2,0){8}{\circle*{.2}}
\multiput(-14,0)(2,0){14}{\circle*{.2}}
\multiput(-15,-1)(2,0){16}{\circle*{.2}}
\multiput(-18,-2)(2,0){19}{\circle*{.2}}
\multiput(-19,-3)(2,0){4}{\circle*{.2}}
\multiput(13,-3)(2,0){4}{\circle*{.2}}
\multiput(-18,-4)(2,0){3}{\circle*{.2}}
\multiput(14,-4)(2,0){3}{\circle*{.2}}
\multiput(-19,-5)(2,0){3}{\circle*{.2}}
\multiput(15,-5)(2,0){3}{\circle*{.2}}
\multiput(-18,-6)(2,0){3}{\circle*{.2}}
\multiput(14,-6)(2,0){3}{\circle*{.2}}
\multiput(-19,-7)(2,0){4}{\circle*{.2}}
\multiput(13,-7)(2,0){4}{\circle*{.2}}
\multiput(-18,-8)(2,0){19}{\circle*{.2}}
\multiput(-15,-9)(2,0){16}{\circle*{.2}}
\multiput(-14,-10)(2,0){14}{\circle*{.2}}
\multiput(-7,-11)(2,0){8}{\circle*{.2}}
\multiput(-2,-12)(2,0){3}{\circle*{.2}}
\put(15,4){\linethickness{1pt}\vector(1,1){20}}
\put(-15,4){\linethickness{1pt}\vector(-1,1){20}}
\put(-40,35){\linethickness{.8pt}\qbezier(20,0)(20,8)(0,8)
\qbezier(-20,0)(-20,8)(0,8) \qbezier(20,0)(20,-8)(0,-8)
\qbezier(-20,0)(-20,-8)(0,-8)}
\put(40,35){\linethickness{.8pt}\qbezier(20,0)(20,8)(0,8)
\qbezier(-20,0)(-20,8)(0,8) \qbezier(20,0)(20,-8)(0,-8)
\qbezier(-20,0)(-20,-8)(0,-8)}
\put(-55,34){SU(3) holonomy} \put(28,34){nearly K\"ahler}
\put(-15,-5){\linethickness{.8pt}\qbezier(35,0)(35,8)(0,8)
\qbezier(-35,0)(-35,8)(0,8) \qbezier(35,0)(35,-8)(0,-8)
\qbezier(-35,0)(-35,-8)(0,-8)}
\put(15,-5){\linethickness{.8pt}\qbezier(35,0)(35,8)(0,8)
\qbezier(-35,0)(-35,8)(0,8) \qbezier(35,0)(35,-8)(0,-8)
\qbezier(-35,0)(-35,-8)(0,-8)}
\put(-38,-6){hypo} \put(-12,-6){double hypo} \put(24,-6){nearly
hypo}
\put(29,14){nearly hypo} \put(29,10){evolution equations}
\put(-62,14){Conti-Salamon} \put(-62,10){evolution equations}
\end{picture}

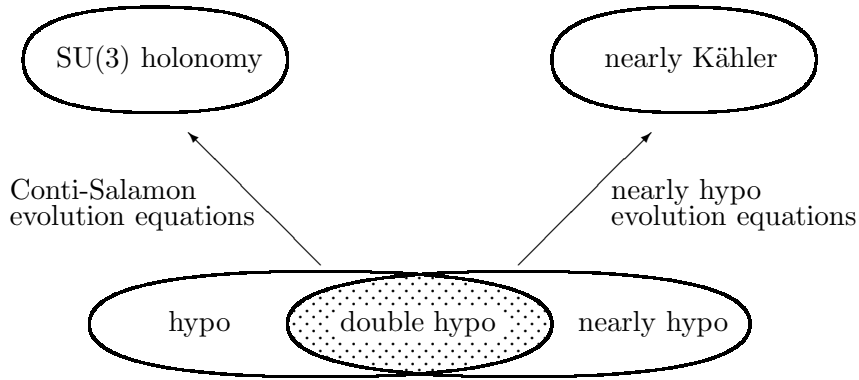
\centerline{{\bf Figure 2}: Special metrics obtained from evolution
of SU(2)-structures}

\bigskip

In Figure $2$, we write $SU(3)$ holonomy for $SU(3)$-structures such
that the holonomy of its metric is contained in $SU(3)$. Moreover,
taking into account $(3.3)$, we must notice that the sin-cone metric of a
Sasaki-Einstein structure on a $5$-manifold $N^5$ defines a nearly
K\"ahler metric on $N^5\times\mathbb R$ and, by $(2.13)$, the cone
metric of a Sasaki-Einstein structure on $N^5$ defines a metric  on
$N^5\times\mathbb R$ whose holonomy is contained in $SU(3)$.

In order to give a characterization of double hypo structures, we
first recall that an $SU(3)$-structure $(F,\Psi_+,\Psi_-)$ on a
$6$-manifold $M^6$ is called half-flat if it satisfies the
conditions
 \begin{equation}\label{hflat}
 dF\wedge F=
 d\Psi_+ =0.
 \end{equation}
These structures become of recent interest mainly because Hitchin
shows in \cite{Hit} that an $SU(3)$-structure on $M^6$ can be lifted
to a $G_2$-holonomy structure on $M^6\times\mathbb R$,  exactly when
the underlying $SU(3)$-structure is half flat.

\begin{thm}\label{dhalf}
Let $(\eta,\omega_i)$ be an $SU(2)$-structure on a $5$-manifold
$N^5$. The following conditions are equivalent:
\begin{itemize}
\item[i).]  $(\eta,\omega_i)$ is a double hypo structure;
\item[ii).] $(\eta,\omega_i)$ satisfies the equations
\begin{equation}\label{dhypo}
d(\eta\wedge\omega_1)=0,\quad d\omega_1=3\eta\wedge\omega_2,\quad
d(\eta\wedge\omega_3)=-2\omega_1\wedge\omega_1, \quad d\omega_3=0.
\end{equation}
\item[iii).] the sin-cone $(N^5\times (0,\pi),F,\Psi_+,\Psi_-)$ with the
$SU(3)$-structure determined by \eqref{newnks} is half-flat.
\end{itemize}
\end{thm}
\begin{proof}
The equivalence of i) and ii) is straightforward consequence from
\eqref{rhypo}, \eqref{nkhypo} and \eqref{dhypo}.

We calculate from \eqref{newnks} that
$$d\Psi_+=\sin^2t\, dt\wedge[\cos t(3\eta\wedge\omega_2-d\omega_1)+\sin t\,
d\omega_3]+\sin^3t\, d(\eta\wedge\omega_2).$$ Consequently,
$d\Psi_+=0 \Longrightarrow d\omega_1=3\eta\wedge\omega_2, \quad
d\omega_3=0$. Using this equivalence, we obtain
$$d(F\wedge F)=2\sin^3t[\cos t(2\omega_1\wedge\omega_1+d(\eta\wedge\omega_3))+\sin t\,
\omega_1\wedge d\eta]\wedge dt.$$ Hence, \eqref{dhypo} are
equivalent to \eqref{hflat}
\end{proof}

\subsection{Double hypo structures on Lie groups}

Next we determine the left-invariant double hypo structures on Lie
groups $G$ satisfying $[\frg,\frg] \not= \frg$, where $\frg$ denotes
the Lie algebra of $G$. In particular we show that solvable Lie
groups cannot admit structures of this type.

\begin{pro}\label{red-equations}
Let $G$ be a Lie group endowed with a left-invariant double hypo
structure $(\eta, \omega_i)$. If the Lie algebra $\frg$ of $G$
satisfies $[\frg,\frg] \not= \frg$, then there is a basis
$e^1,\ldots,e^5$ for ${\frak g}^*$ and a real number $\mu$ such that
\begin{equation}\label{structure}
\eta=e^5,\quad\quad \omega_1= e^{12} + e^{34},\quad\quad \omega_2=
e^{13} + e^{42},\quad\quad \omega_3= e^{14} + e^{23},
\end{equation}
and
\begin{equation}\label{differentials}
\left\{
\begin{array}{rl}
d e^1= \!\!\!&\!\! 0,\\[3pt]
d e^2= \!\!\!&\!\! \mu e^{34} -3 e^{35},\\[4pt]
d e^3= \!\!\!&\!\! - \mu e^{24} + 3 e^{25},\\[4pt]
d e^4= \!\!\!&\!\! \mu e^{14},\\[4pt]
d e^5 = \!\!\!&\!\! -4 e^{23} + {\mu^2\over 3} (e^{14} - e^{23}).
\end{array}
\right.
\end{equation}
\end{pro}

\begin{proof}
Let $V$ be the subspace of $\frg^*$ orthogonal to $\eta$, and let
$\alpha\in \frg^*$ be nonzero and closed. Thus, $\alpha = \beta +
\rho\, \eta$, where $\beta\in V$ and $\rho\in \mathbb{R}$. Now,
$d\alpha=0$ is equivalent to $d\beta = - \rho\, d\eta$. Therefore,
$\gamma={1\over ||\beta||}\, \beta$ is a unit element in $V=\langle
\eta \rangle^\perp$ satisfying
$$
d\gamma = \lambda\, d\eta,
$$
with $\lambda =-\rho/||\beta||$. By~\cite[Corollary 3]{ConS} there
is a basis $e^1,\ldots,e^5$ for ${\frak g}^*$ satisfying
(\ref{structure}) with $e^1=\gamma$. In terms of this basis the
differentials of $e^1,\ldots,e^5$ are given by
\begin{equation}\label{gen-differentials}
d e^i= \sum_{1\leq j<k\leq 5} c^i_{jk}\, e^{jk}, \quad\quad 1\leq
i\leq 5,
\end{equation}
where $c^1_{jk}=\lambda\, c^5_{jk}$ for all $j,k$. The 41 remaining
coefficients $\lambda, c^2_{jk},\ldots,c^5_{jk}$ must satisfy the
Jacobi identity $d(de^i)=0$, $1\leq i\leq 5$, and the double hypo
conditions~(\ref{dhypo}).

First, a direct calculation shows that
$$
\begin{array}{rl}
 d\omega_1 =  \!\!\!&\!\! d e^{12}+ d e^{34}= -(\lambda c^5_{13} + c^2_{23}+c^4_{12}) e^{123} -
(\lambda c^5_{14} + c^2_{24}-c^3_{12}) e^{124}- (\lambda c^5_{15} +
c^2_{25})
e^{125}\\[5pt]
&  - (c^2_{34}- c^3_{13}-c^4_{14}) e^{134}-(c^2_{35}-c^4_{15})
e^{135}
-(c^2_{45}+c^3_{15}) e^{145} +(\lambda c^5_{34} + c^3_{23}+c^4_{24}) e^{234}\\[5pt]
& +(\lambda c^5_{35} + c^4_{25}) e^{235} +(\lambda c^5_{45} -
c^3_{25}) e^{245} - (c^3_{35} + c^4_{45}) e^{345}.
\end{array}
$$
Since $3\eta\wedge\omega_2 = 3 e^{135} - 3 e^{245}$, we have that
$d\omega_1=3\eta\wedge\omega_2$ if and only if the coefficients
$\lambda, c^i_{jk}$ satisfy the following relations:
\begin{equation}\label{coeffs-1}
\begin{array}{rll}
& c^2_{25}= -\lambda\, c^5_{15},\ \, c^3_{12}= \lambda\, c^5_{14} +
c^2_{24},\ \, c^3_{15}= - c^2_{45},\ \, c^3_{25} = 3 + \lambda\,
c^5_{45},\ \, c^4_{12}= -\lambda\,
c^5_{13} -c^2_{23},\\[5pt]
& c^4_{14}= c^2_{34} -c^3_{13},\ \, c^4_{15}= 3 + c^2_{35},\ \,
c^4_{24}= -\lambda\, c^5_{34}- c^3_{23},\ \, c^4_{25}= -\lambda\,
c^5_{35},\ \, c^4_{45}= -c^3_{35}.
\end{array}
\end{equation}

On the other hand, since
$$
\begin{array}{rl}
 d(\eta\wedge\omega_3) =  \!\!\!&\!\! de^{145}+de^{235}= (c^5_{14} + c^5_{23}) e^{1234}
 + (c^5_{15} + c^2_{12} + c^3_{13} - c^4_{23}) e^{1235}\\[5pt]
& + (\lambda\, c^5_{12} - c^5_{25} + \lambda\, c^5_{34} + c^3_{14} +
c^3_{23})
e^{1245} + (\lambda\, c^5_{13} - c^5_{35} - c^2_{14}- c^4_{34}) e^{1345}\\[5pt]
&  + (\lambda\, c^5_{23} + c^5_{45} - c^2_{24}- c^3_{34}) e^{2345}
\end{array}
$$
and $\omega_1\wedge\omega_1= 2 e^{1234}$, we conclude that
$d(\eta\wedge\omega_3)=-2\omega_1\wedge\omega_1$ if and only if
\begin{equation}\label{coeffs-2}
\begin{array}{rl}
& c^3_{23}= -\lambda\, c^5_{12} + c^5_{25}-\lambda\, c^5_{34} -
c^3_{14},\quad c^3_{34}= \lambda\, c^5_{23}
+c^5_{45} - c^2_{24},\quad  c^4_{23}= c^5_{15} +c^2_{12} +c^3_{13},\\[5pt]
& c^4_{34}= \lambda\, c^5_{13}-c^5_{35} - c^2_{14},\quad c^5_{23}=
-4-c^5_{14}.
\end{array}
\end{equation}

A direct calculation shows that
$$
\begin{array}{rl}
 d\omega_3 = \!\!\!&\!\! de^{14}+de^{23}= - c^5_{15} e^{123}
 + c^5_{25} e^{124} + (\lambda\, c^5_{35} - c^2_{45}) e^{125}
 + c^5_{35} e^{134}- (c^2_{15} + c^4_{35}) e^{135}\\[5pt]
& -(\lambda\, c^5_{15}-c^3_{35}) e^{145} - c^5_{45} e^{234} +
(\lambda\, c^5_{15} - c^3_{35}) e^{235} - (\lambda\, c^5_{25} +
c^3_{45}) e^{245} -(\lambda\, c^5_{35}-c^2_{45}) e^{345},
\end{array}
$$
which implies that $\omega_3$ is closed if and only if
\begin{equation}\label{coeffs-21}
c^2_{45}=c^3_{35}=c^3_{45}=c^5_{15}=c^5_{25}=
c^5_{35}=c^5_{45}=0,\quad\quad\quad c^4_{35}= -c^2_{15}.
\end{equation}

Moreover, the 3-form $\eta\wedge \omega_1$ is closed if and only if
the additional relation
\begin{equation}\label{coeffs-22}
c^5_{34}=-c^5_{12}
\end{equation}
is satisfied.

Notice also that
$0=d^2(\omega_1)=3d(\eta\wedge\omega_2)=3(d\eta\wedge\omega_2
-\eta\wedge d\omega_2)$ implies $d\eta\wedge\omega_2 =\eta\wedge
d\omega_2$, which is equivalent to the conditions
\begin{equation}\label{coeffs-3}
c^2_{23}= - c^2_{14},\ \ c^2_{24}= -4\lambda +
 c^2_{13},\ \  c^3_{24}= -c^2_{12} -c^2_{34}
 +c^3_{13},\ \ c^4_{13}= \lambda\, c^5_{12} + c^3_{14},\ \ c^5_{24}=
c^5_{13}.
\end{equation}
Now, it is easy to see that the coefficient of $e^{245}$ in the
3-form $d(de^5)$ vanishes if and only if
\begin{equation}\label{coeffs-31}
c^5_{12}= 0.
\end{equation}

>From (\ref{coeffs-1})--(\ref{coeffs-31}) we get that the structure
equations (\ref{gen-differentials}) reduce to
\begin{equation}\label{differentials-bis}
\left\{\!\!
\begin{array}{rl}
d e^1= \!\!\!&\!\! \lambda\, d e^5,\\[5pt]
d e^2= \!\!\!&\!\! c^2_{12} e^{12} + c^2_{13} e^{13} + c^2_{14}
e^{14} + c^2_{15} e^{15} - c^2_{14} e^{23} - (4\lambda - c^2_{13})
e^{24} + c^2_{34}
e^{34}+ c^2_{35} e^{35},\\[5pt]
d e^3= \!\!\!&\!\! -(4\lambda - \lambda\, c^5_{14} - c^2_{13})
e^{12} + c^3_{13} e^{13} + c^3_{14} e^{14} - c^3_{14} e^{23}   -
(c^2_{12} + c^2_{34} - c^3_{13}) e^{24}\\[4pt]
& + 3 e^{25} - (\lambda\, c^5_{14} + c^2_{13}) e^{34},\\[5pt]
d e^4= \!\!\!&\!\! -(\lambda\, c^5_{13} -c^2_{14}) e^{12} + c^3_{14}
e^{13} + (c^2_{34} - c^3_{13}) e^{14}
+ (3 + c^2_{35}) e^{15} + (c^2_{12} + c^3_{13}) e^{23}\\[4pt]
& +  c^3_{14} e^{24} + (\lambda\, c^5_{13} -
c^2_{14}) e^{34} - c^2_{15} e^{35},\\[5pt]
d e^5 = \!\!\!&\!\!  c^5_{13} e^{13} + c^5_{14} e^{14} - (4+
c^5_{14}) e^{23} + c^5_{13} e^{24},
\end{array}
\right.
\end{equation}
where the 11 coefficients $\lambda, c^2_{12}, c^2_{13}, c^2_{14},
c^2_{15}, c^2_{34}, c^2_{35}, c^3_{13}, c^3_{14}, c^5_{13}$ and
$c^5_{14}$  must satisfy the Jacobi identity $d(de^i)=0$, for $1\leq
i \leq 5$.

For the rest of the proof we follow a decision tree depending on the
nullity of the coefficients to conclude that the Jacobi identity is
satisfied if and only if $c^5_{14}=(c^2_{34})^2/3$ and the remaining
coefficients vanish. The proof of this fact is rather long but
straighforward, so we omit details.
\end{proof}

It is easy to see that the Lie group determined by
(\ref{differentials}) is isomorphic to $SU(2) \times A^2$ for $\mu=0$
and $SU(2)\times A\!f\!f(\RR)$ for $\mu\not=0$, where $A^2$ denotes
a 2-dimensional abelian Lie group and $A\!f\!f(\RR)$ is the group of
affine transformations of $\RR$. As a consequence of
Proposition~\ref{red-equations}, the Lie group $SU(2)\times A^2$ has
a unique (up to equivalence) left-invariant double hypo structure.
Moreover:

\begin{co}\label{family1}
Let $\{X_1,X_2,X_3\}$ and $\{Y_1,Y_2\}$ be the standard basis of
left-invariant vector fields on 
$SU(2)$
and $A^2$, respectively,
that is,
$$
[X_1,X_2]=X_3,\quad [X_2,X_3]=X_1,\quad [X_3,X_1]=X_2,\quad
[Y_1,Y_2]=0,
$$
and let us denote by $\{\alpha^i,\beta^j\}$ the dual basis of
$\{X_i,Y_j\}$. For each $\rho\in\RR$, the 
$SU(2)$-structure on the Lie group 
$SU(2)\times A^2$ given by
$$
\begin{array}{rlll}
\!\!\!&\!\!  \eta= {1\over 3}\alpha^1,&\quad \!\!\!&\!\! \omega_1=
{1\over 2\sqrt{3}}(-\alpha^2\beta^1+{\rho\over 6\sqrt{3}}
\alpha^2\alpha^3+{1\over 3}\alpha^3\beta^2),\\[9pt]
\!\!\!&\!\!  \omega_2= -{1\over 2\sqrt{3}}(\alpha^3\beta^1+{1\over
3}\alpha^2\beta^2),&\quad \!\!\!&\!\! \omega_3= {\rho\over
6\sqrt{3}}\alpha^2\beta^1+ {1\over 3} \beta^1\beta^2+ {1\over
12}\alpha^2\alpha^3,
\end{array}
$$
is nearly hypo, and it is double hypo if and only if $\rho=0$.
\end{co}

\begin{proof}
In terms of the new basis $\{e^i\}$ defined by
$$\alpha^1=3\,e^5,\quad \alpha^2=2\sqrt{3}\, e^2,\quad
\alpha^3=2\sqrt{3}\, e^3,\quad \beta^1=e^1,\quad \beta^2 =\rho\,
e^2+3\,e^4,$$ the structure equations of the Lie group are
$$
de^1=0,\quad de^2=-3\,e^{35},\quad de^3=3\,e^{25},\quad
de^4=\rho\,e^{35},\quad de^5=-4\,e^{23},
$$
and the SU(2)-structure is given by (\ref{structure}). Therefore,
when $\rho=0$ the structure is precisely the one given in
Proposition~\ref{red-equations} for $\mu=0$. Finally, it is easy to
check that for each $\rho\not=0$ the structure is nearly hypo but
the form $\omega_3$ is not closed.
\end{proof}

>From our discussion above and Remark~\ref{answer-Q1} it follows in
particular that $S^3\times \TT^2$ is a real analytic manifold having
an analytic double hypo structure, therefore:

\begin{co}
The double hypo structure on $S^3\times \TT^2$ given by
$$
\eta= {1\over 3}\alpha^1,\quad \omega_1= -{1\over
2\sqrt{3}}(\alpha^2\beta^1-{1\over 3}\alpha^3\beta^2),\quad
 \omega_2= -{1\over 2\sqrt{3}}(\alpha^3\beta^1+{1\over
3}\alpha^2\beta^2),\quad  \omega_3= {1\over 3} \beta^1\beta^2+
{1\over 12}\alpha^2\alpha^3,
$$
can be lifted both to a nearly K\"ahler structure and to a
Calabi-Yau structure.
\end{co}

>From Proposition~\ref{red-equations} we get that solvable Lie groups
do not admit left-invariant double hypo structures. Since there
exist nilpotent Lie groups having left-invariant hypo
structures~\cite{ConS}, the class of manifolds with double hypo
structures is a proper subclass of that consisting of hypo
manifolds. Moreover, 
Corollary~\ref{family1} shows the existence of
nearly hypo structures which are not double hypo.

\begin{co}\label{family2}
The Lie group
 $SU(2)\times A\!f\!f(\RR)$ admits a
$1$-parametric family of left-invariant double hypo structures. More
precisely, if $\{X_1,X_2,X_3\}$ and $\{Z_1,Z_2\}$ are the standard
basis of left-invariant vector fields on 
$SU(2)$ and $A\!f\!f(\RR)$, respectively, that is,
$$
[X_1,X_2]=X_3,\quad [X_2,X_3]=X_1,\quad [X_3,X_1]=X_2,\quad
[Z_1,Z_2]=Z_2,
$$
then (up to equivalence) any left-invariant double hypo structure
$(\eta,\omega_i)$ on 
$SU(2)\times A\!f\!f(\RR)$ belongs to the
family
$$
\begin{array}{rlll}
\!\!\!&\!\!  \eta= {1\over 3}(\alpha^1+\mu\,\gamma^2),&\quad
\!\!\!&\!\! \omega_1= {1\over \mu(\mu^2+12)^{1\over
2}}(\alpha^2\gamma^1+\mu\,
\alpha^3\gamma^2),\\[10pt]
\!\!\!&\!\!  \omega_2= {1\over \mu(\mu^2+12)^{1\over
2}}(\alpha^3\gamma^1-\mu\, \alpha^2\gamma^2),&\quad  \!\!\!&\!\!
\omega_3= -{1\over \mu}\gamma^1\gamma^2+ {1\over \mu^2+12}
\alpha^2\alpha^3,
\end{array}
$$
for some $\mu\in \RR-\{0\}$, where $\{\alpha^i,\gamma^j\}$ denotes
the dual basis of $\{X_i,Z_j\}$.
\end{co}

\begin{proof}
It follows directly from Proposition~\ref{red-equations} taking
$$\alpha^1=-\mu e^4+3e^5,\quad \alpha^2=(\mu^2+12)^{1\over 2}\, e^2,\quad
\alpha^3=(\mu^2+12)^{1\over 2}\, e^3,\quad \gamma^1=-\mu e^1,\quad
\gamma^2 =e^4.$$
\end{proof}

We finish this section by showing that the double hypo structures of
Proposition~\ref{red-equations} can be deformed into hypo
structures. In fact, it is easy to see that for each $r\in
\RR-\{0\}$ and $\mu,\tau\in\RR$, the Lie group $G$ determined by the
equations
\begin{equation}\label{deform}
\left\{
\begin{array}{rl}
d e^1= \!\!\!&\!\! 0,\\[3pt]
d e^2= \!\!\!&\!\! \mu e^{34} +r e^{35},\\[4pt]
d e^3= \!\!\!&\!\! - \mu e^{24} - r e^{25},\\[4pt]
d e^4= \!\!\!&\!\! \mu e^{14},\\[4pt]
d e^5 = \!\!\!&\!\! (\tau-{\mu^2\over r}) e^{23} - {\mu^2\over r}
(e^{14} - e^{23}),
\end{array}
\right.
\end{equation}
is isomorphic to the product $H\times K$, where $H=A^2$ for $\mu=0$
and $H=A\!f\!f(\RR)$ for $\mu\not=0$, and $K=$ SU(2) if $r\tau>0$,
$K=$ SL(2, $\RR$) if $r\tau<0$ and $H=$ E(2) if $\tau=0$, E(2) being
the group of rigid motions of Euclidean 2-space. Moreover, a direct
calculation shows that the SU(2)-structure given
by~(\ref{structure}) is always hypo, and it is double hypo if and
only if $r=-3$ and $\tau=-4-{\mu^2\over 3}$.

\section{Nearly half flat structures on $6$-manifolds}\label{nearlyhalf}

In this section we generalize the construction of nearly parallel
$G_2$-structures on $M^6\times\mathbb R$ induced from a nearly
K\"ahler structure on $M^6$ described in \cite{BM}. For general
results on $G_2$-manifolds, see \cite{FG}.

Let $(F,\Psi_+,\Psi_-)$ be an $SU(3)$-structure on a $6$-manifold
$M^6$. We consider the $G_2$-structure $\phi$ on $M^6\times\mathbb
R$ defined by the $3$-form $\phi$ given by
  \begin{equation}\label{g21}
  \phi=F\wedge dq - \Psi_-,
  \end{equation}
where $dq$ is the standard $1$-form on $\RR$. We also have a
$4$-form
  \begin{equation}\label{g22}
  \Hodge_7\, \phi=\frac12F\wedge F + \Psi_+\wedge dq,
  \end{equation}
where $\Hodge_7$ denotes the Hodge star operator on
$M^6\times\mathbb R$.

Vice versa, let $f: M^6\longrightarrow P^7$ be a hypersurface of a
$G_2$-manifold $(P^7,\phi)$ and denote by $\mathbb N$ the unit
normal. Then the $G_2$-structure $\phi$ induces an $SU(3)$-structure
$(F,\Psi_+,\Psi_-)$ on $M^6$ defined by the equalities
  \begin{equation}\label{hypg21}
  F=\mathbb N\lrcorner \phi,\quad \Psi_+=-\mathbb
  N\lrcorner\Hodge\phi, \quad \Psi_-=-f^*\phi.
  \end{equation}
The types of the induced $U(3)$-structures are investigated in
\cite{Cal,Gr,Br1} while the types of the induced $SU(3)$-structures
are studied recently in \cite{Cab1}.

We recall that a $G_2$-structure is called {\em nearly parallel} if
  \begin{equation}\label{g2np}
  d\phi=4\Hodge\phi.
  \end{equation}
It is well known that nearly parallel $G_2$-structures are Einstein
with positive scalar curvature $s=54\cdot 7\cdot 16=6048$.

Hitchin shows in \cite{Hit} that an $SU(3)$-structure on $M^6$ can
be lifted to a {\em parallel} $G_2$-structure  \eqref{g21} on
$M^6\times\mathbb R$, i.e. \cite{FG}, a $G_2$-structure satisfying
$d\phi=d\Hodge\phi=0$ (or, equivalently,  $M^6\times\mathbb R$ has a
metric whose holonomy is contained in $G_2$), exactly when the
underlying $SU(3)$-structure is half flat (note that the half-flat
condition compatible with \eqref{g21} reads $dF\wedge F=d\Psi_-=0$).
Thus, any double hypo structure on a 5-manifold could produce a
$G_2$-holonomy metric  by solving Hitchin's flow equations
(compatible with \eqref{g21})
$$\partial_q\Psi_-=-dF, \quad d\Psi_+=-\frac12\partial_q(F\wedge F)$$
since its sin-cone is half-flat due to Theorem \ref{dhalf}.

Next, we search for sufficient conditions imposed  on an
$SU(3)$-structure $(F,\Psi_+,\Psi_-)$ which imply that the
$G_2$-structure  on $M^6\times \mathbb R$ determined by \eqref{g21}
is nearly parallel, i.e. it satisfies \eqref{g2np}.

\begin{defn}
We call an $SU(3)$-structure $(F,\Psi_+,\Psi_-)$ on a $6$-manifold
$M^6$ \emph{nearly half flat} if it satisfies the equation
  \begin{gather}\label{g2hypo}
  d\Psi_-=-2F\wedge F.
  \end{gather}
In particular, any nearly K\"ahler $6$-manifold carries a nearly
half flat structure.
\end{defn}

The following diagram represents the relations among
$SU(3)$-structures on $6$-manifolds:

\setlength{\unitlength}{2.5pt}
\begin{picture}(120,70)(-82,-18)
\put(-15,3){\linethickness{1pt}\line(1,0){31}}
\put(-15,3){\linethickness{1pt}\line(0,1){10}}
\put(-15,13){\linethickness{1pt}\line(1,0){31}}
\put(16,3){\linethickness{1pt}\line(0,1){10}}
\put(-12.5,7){nearly K\"ahler}
\put(-25,-15){\linethickness{1pt}\line(1,0){90}}
\put(-25,-15){\linethickness{1pt}\line(0,1){35}}
\put(-25,20){\linethickness{1pt}\line(1,0){90}}
\put(65,-15){\linethickness{1pt}\line(0,1){35}}
\put(-65,-4){\linethickness{1pt}\line(1,0){90}}
\put(-65,-4){\linethickness{1pt}\line(0,1){35}}
\put(-65,31){\linethickness{1pt}\line(1,0){90}}
\put(25,-4){\linethickness{1pt}\line(0,1){35}}
\put(-53,13){half flat}
\put(31,1){nearly half flat}
\put(-75,-20){\linethickness{1pt}\line(1,0){150}}
\put(-75,-20){\linethickness{1pt}\line(0,1){68}}
\put(-75,48){\linethickness{1pt}\line(1,0){150}}
\put(75,-20){\linethickness{1pt}\line(0,1){68}}
\put(-30,38){6-manifolds with $SU(3)$-structure}
\end{picture}

\bigskip

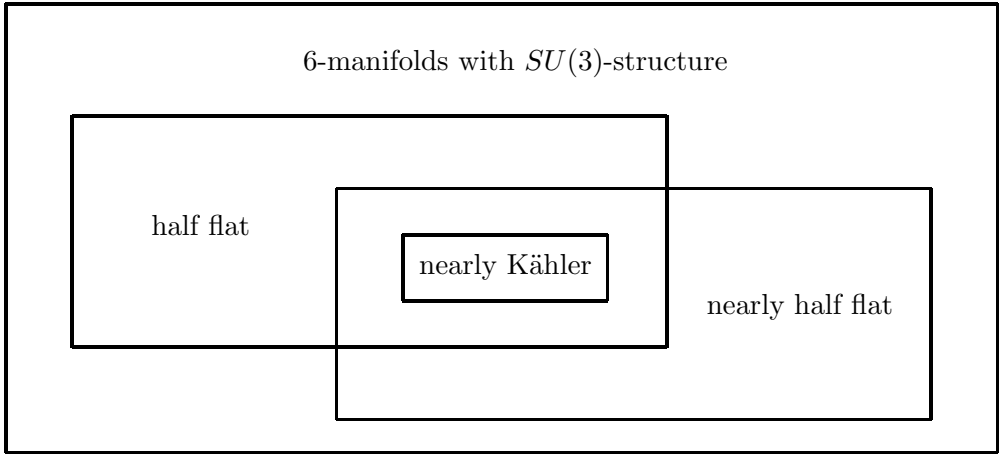
\centerline{{\bf Figure 3}: Classes of $SU(3)$-structures}

\medskip
Consider $SU(3)$-structures $(F(q),\Psi_+(q),\Psi_-(q))$ on $M^6$
depending on a real parameter $q\in\mathbb R$ and the corresponding
$G_2$-structure $\phi(q)$ on $M^6\times\mathbb R$. We have

\begin{pro}\label{g2evpro}
An $SU(3)$-structure $(F,\Psi_+,\Psi_-)$ on $M^6$ can be lifted to a
nearly parallel $G_2$-structure $\phi(q)$ on $M^6\times\mathbb R$
defined by \eqref{g21} if and only if it is a nearly half flat
structure and the following \emph{evolution nearly half flat
equations} hold
  \begin{equation}\label{evolug2}
  \begin{array}{c}
  \partial_q\Psi_- =4\Psi_+-dF,\\
  d\Psi_+=-\frac12\partial_q(F\wedge F).
  \end{array}
  \end{equation}
\end{pro}
\begin{proof}
Take the exterior derivative in \eqref{g21} and use \eqref{g22} to
get that the equation \eqref{g2np} holds precisely when
\eqref{g2hypo} and \eqref{evolug2} are fulfilled. Moreover,
\eqref{evolug2} imply that \eqref{g2hypo} holds for any time $q$ due
to the equality $\partial_q(d\Psi_-+2F\wedge
F)=4d\Psi_++2\partial_q(F\wedge F)$.
\end{proof}

As a consequence of the above considerations, we can recover one of
the main results in~\cite{BM}.

\begin{thm}\cite{BM}\label{newg2}
Let $(M^6,F,\Psi_+,\Psi_-)$ be a nearly K\"ahler structure. Then the
$G_2$-structure $\phi$ on $M^6\times\mathbb R$ defined for $0\le
q\le \pi$ by
  \begin{gather}\label{oldg2}
  \phi=\sin^2q\, F\wedge dq -\sin^3q\left(-\cos q\, \Psi_+ + \sin
  q\,\Psi_-\right)
  \end{gather}
is a nearly parallel $G_2$-structure on $M^6\times\mathbb R$
generating the well known Einstein metric
$$g_7=dq^2+\sin^2q\, g_6,$$
where $g_6$ is the nearly K\"ahler metric on $M^6$.

If $(M^6,F,\Psi_+,\Psi_-)$ is compact then $(M^6\times
[0,\pi],\phi)$ is a compact nearly parallel $G_2$-manifold with two
conical singularities at $q=0$ and  $q=\pi$.
\end{thm}

\begin{proof}
Consider the  $SU(3)$-structure $(F(q),\Psi_+(q),\Psi_-(q))$
depending on a real parameter $q$:
  \begin{equation}\label{evog2}
  \begin{aligned}
  &F(q)=\sin^2q\, F \\
  &\Psi_+(q)=\sin^3q\left(\sin q\, \Psi_+ + \cos q\, \Psi_-\right),\\
  &\Psi_-(q)=\sin^3q\left(-\cos q\, \Psi_+ + \sin q\,\Psi_-\right).
  \end{aligned}
  \end{equation}
Applying \eqref{nkdef}, we see that the structure defined by
\eqref{evog2} satisfies the nearly half flat  conditions
\eqref{g2hypo} as well as the evolution nearly parallel equation
\eqref{evolug2}. Consequently, the structure \eqref{oldg2} satisfies
\eqref{g2np} and therefore it is a nearly parallel $G_2$-structure
on $M^6\times\mathbb R$.
\end{proof}

As a consequence of the proof of Theorem \ref{newg2}, we obtain
\begin{co}\label{conewg2}
An $SU(3)$-manifold  $(N^6,F,\Psi_+,\Psi_-)$ is nearly K\"ahler if
and only if  the sin-cone $(N^6\times\mathbb R,\phi)$ with the
$G_2$-structure  defined by \eqref{oldg2} is a nearly parallel
$G_2$-manifold for any $0< t<\pi$.
\end{co}
\begin{proof}
The equations \eqref{oldg2} imply
$$d\phi=\sin^3q\cos q\, d\Psi_+-\sin^4q\,d\Psi_-+ [\sin^2q\,dF-
(3\sin^2q\cos^2q-\sin^4q)\Psi_++4\sin^3q\cos q\,\Psi_-]\wedge dq.
$$
Consequently, $d\phi=4*\phi \Leftrightarrow
d\omega_1=3\eta\wedge\omega_2, \quad d\eta=-2\omega_3.$ Using this
equivalence, we obtain
$$d\phi-4*\phi =\sin^3q[\cos q\,d\Psi_+-\sin q(d\Psi_-+2F\wedge F)] +\sin^2q(dF-3\Psi_+)\wedge dq.
$$
Hence, $d\phi=4*\phi  \Leftrightarrow dF=3\Psi_+, \quad
d\Psi_-=-2F\wedge F$. Thus, \eqref{nkdef} are equivalent to
\eqref{g2np}  and the proof is complete.
\end{proof}

More generally we have
\begin{pro}\label{hypg2new}
Let $f: M^6\longrightarrow P^7$ be an immersion of an oriented
$6$-manifold into a $7$-manifold with a nearly parallel
$G_2$-structure. Then the $SU(3)$-structure induced on $M^6$ is a
nearly half flat  $SU(3)$-structure.
\end{pro}

\begin{proof}
Since $f^*$ commutes with $d$, the equalities \eqref{hypg21}
substituted into \eqref{g2np} yield \eqref{g2hypo}.
\end{proof}

{\bf Question 2.} Does the converse of Proposition~\ref{hypg2new}
hold? i.e. is it true that any (real analytic) nearly half flat
structure on $M^6$ can be lifted to a nearly parallel
$G_2$-structure on $M^6\times\RR$? This is equivalent to prove the
existence of a solution of the evolution nearly half flat equation
\eqref{evolug2}. \footnote{Recently we learned that Stock proves in
Theorem $2.5$ of \cite{Stock} that Question $2$ has an affirmative answer
for nearly half flat structures on closed $6$-manifolds $M^6$, i.e.
they can be lifted to a nearly parallel $G_2$-structure on
$M^6\times I$, for a sufficiently small interval $I$.}

Notice that nearly K\"ahler structures can be lifted, on one hand,
to a metric with holonomy contained in $G_2$ (that is, to a parallel
$G_2$-structure) due to Hitchin result \cite{Hit} and, on the other
hand, taking account Corollary \ref{conewg2},
to a nearly parallel $G_2$-structure, providing a relation between
these special classes on 7-dimensional manifolds:

\setlength{\unitlength}{2.5pt}
\begin{picture}(120,70)(-82,-18)
\put(15,4){\linethickness{1pt}\vector(1,1){20}}
\put(-15,4){\linethickness{1pt}\vector(-1,1){20}}
\put(-40,35){\linethickness{.8pt}\qbezier(20,0)(20,8)(0,8)
\qbezier(-20,0)(-20,8)(0,8) \qbezier(20,0)(20,-8)(0,-8)
\qbezier(-20,0)(-20,-8)(0,-8)}
\put(40,35){\linethickness{.8pt}\qbezier(20,0)(20,8)(0,8)
\qbezier(-20,0)(-20,8)(0,8) \qbezier(20,0)(20,-8)(0,-8)
\qbezier(-20,0)(-20,-8)(0,-8)}
\put(-54,34){$G_2$ holonomy} \put(23,34){nearly parallel $G_2$}
\put(-15,-5){\linethickness{.8pt}\qbezier(35,0)(35,8)(0,8)
\qbezier(-35,0)(-35,8)(0,8) \qbezier(35,0)(35,-8)(0,-8)
\qbezier(-35,0)(-35,-8)(0,-8)}
\put(15,-5){\linethickness{.8pt}\qbezier(35,0)(35,8)(0,8)
\qbezier(-35,0)(-35,8)(0,8) \qbezier(35,0)(35,-8)(0,-8)
\qbezier(-35,0)(-35,-8)(0,-8)}
\put(0,-5){\linethickness{.8pt}\qbezier(15,0)(15,4)(0,4)
\qbezier(-15,0)(-15,4)(0,4) \qbezier(15,0)(15,-4)(0,-4)
\qbezier(-15,0)(-15,-4)(0,-4)}
\put(-40,-6){half flat} \put(-13,-6){nearly K\"ahler}
\put(21,-6){nearly half flat}
\put(29,14){nearly half flat} \put(29,10){evolution equations}
\put(-62,14){Hitchin evolution} \put(-62,10){equations}
\end{picture}

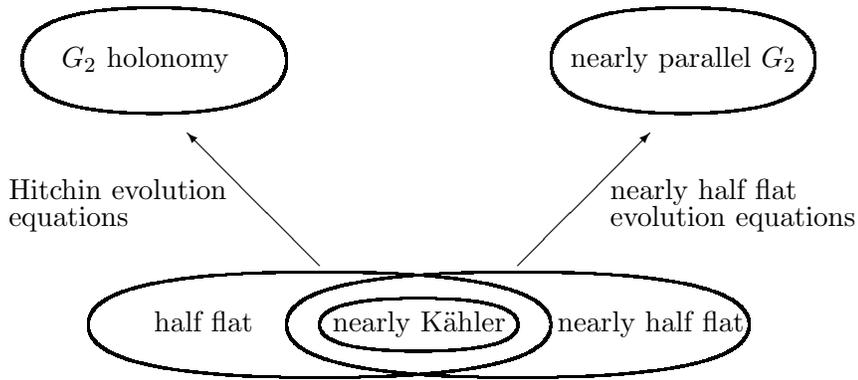
\centerline{{\bf Figure 4}: Special metrics obtained from evolution
of SU(3)-structures}

\section{Examples}\label{examples}

For $N^5=S^5\subset S^6$ and for $N^5=S^2 \times S^3\subset S^3
\times S^3$, we give an explicit description of the Sasaki-Einstein
hypo $SU(2)$-structure on $N^5$ which generates a new nearly
K\"ahler structure with two conical singularities on $S^2 \times
S^3\times [0,\pi]$ as well as a nearly parallel $G_2$-structure on
$N^5\times [0,\pi]\times [0,\pi]$ according to Theorem~\ref{newnk}
and Theorem~\ref{newg2}.  We also apply our results to the new
compact Sasaki-Einstein manifolds $Y^ {p,q}$, which are
diffeomorphic to $S^2\times S^3$ and were constructed recently in
\cite{GMS}, to obtain a new nearly K\"ahler structure with two
conical singularities on $Y^{p,q}\times \RR$ and a nearly parallel
$G_2$-structure on $Y^{p,q}\times \RR^2$.

Finally, we give an example of an analytic double hypo structure and
a solution to the Conti-Salamon hypo evolution equations \eqref{vol}
as well as a solution to the nearly hypo evolution equations
\eqref{evolunk} which is an $SU(2)$-structure only in the beginning
for $t=0$. This shows a difference between Hitchin theorem
\cite{Hit} which says that any solution to the Hitchin flow
equations starting with a half-flat $SU(3)$-structure is
automatically a half-flat $SU(3)$-structure for all $t$.

\subsection{The Nearly K\"ahler structure on $S^5\times\RR$}
We begin with an explicit description of

\subsubsection{The standard $SU(3)$-structure on $S^6$}

Using the stereographic projection of $S^6-\{p\}$  on $\RR^6$ from
the point $p=(0, \cdots,0,1) \in \RR^7$, one can check that a basis
for the vector fields on $S^6-\{p\}$ consists of $\{E_i; 1\leq i\leq
6\}$ with
  \begin{gather*}
  (E_1)_x = (1-x_7-x_1^2, -x_1 x_2, -x_1 x_3, -x_1 x_4, -x_1 x_5, -x_1 x_6, x_1(1-x_7)),\\
  (E_2)_x = (-x_1 x_2, 1-x_7-x_2^2, -x_2 x_3, -x_2 x_4, -x_2 x_5, -x_2 x_6, x_2(1-x_7)),\\
  (E_3)_x = (-x_1 x_3, -x_2 x_3, 1-x_7-x_3^2, -x_3 x_4, -x_3 x_5, -x_3 x_6,  x_3(1-x_7)),\\
  (E_4)_x = (-x_1 x_4, -x_2 x_4, -x_3 x_4, 1-x_7-x_4^2, -x_4 x_5, -x_4 x_6,  x_4(1-x_7)),\\
  (E_5)_x = (-x_1 x_5, -x_2 x_5, -x_3 x_5, -x_4 x_5,  1-x_7-x_5^2, -x_5 x_6, x_5(1-x_7)),\\
  (E_6)_x = (-x_1 x_6, -x_2 x_6, -x_3 x_6, -x_4 x_6, -x_5 x_6, 1-x_7-x_6^2,  x_6(1-x_7)),
  \end{gather*}
for any arbitrary point $x \in S^6-\{p\}$. (Notice that this basis
is orthogonal and $||E_i||^2 = (1-x_7)^2$.) The basis $\{\alpha_i;
1\leq i\leq 6\}$ for the $1$-forms on $S^6-\{p\}$ dual to $\{E_i;
1\leq i\leq 6\}$ is given by
  $$
  \alpha_i= \frac{1}{1-x_7} dx_i +\frac{x_i}{(1-x_7)^2} dx_7.
  $$
{}From now on, we write $x_{ij}=x_{i}x_{j}$,
$x_{ijk}=x_{i}x_{j}x_{k}$, $dx_{ij}=dx_{i} \wedge dx_{j}$, and so
forth. We will need also the expressions of $\alpha_{ij}$ and
$\alpha_{ijk}$ in terms of $dx_{ij}$ and $dx_{ijk}$, respectively;
  $$
  \alpha_{ij} = \frac{1}{(1-x_7)^2} dx_{ij} + \frac{1}{(1-x_7)^3}
  \Big(x_jdx_{i7} - x_idx_{j7}\Big),
  $$
for $1\leq i<j \leq 6$, and
  $$
  \alpha_{ijk} = \frac{1}{(1-x_7)^3} dx_{ijk} + \frac{1}{(1-x_7)^4}
  \Big(x_i dx_{jk7}-x_j dx_{ik7} + x_k dx_{ij7}\Big),
  $$
for $1\leq i < j < k \leq 6$. Let $U=\sum_{i=1}^7 x_i
\frac{\partial}{\partial x_i}$ be the unit normal vector field to
$S^6-\{p\}$. We identify $\RR^7$ with the imaginary part of the
space of Cayley numbers, and define a vector cross product $x \times
y$, where $x$, $y \in \RR^7$, by the imaginary part of the Cayley
number $xy$. Then, the  standard almost complex  structure on $S^6$
is defined by $J(X)=U\times X$ for any vector field X on $S^6$.  A
simple calculation shows that
  \begin{align*}
 (JE_1)_x = \, & (-x_{16}, x_{15} + x_3(1-x_7), x_{14} - x_2(1-x_7), -x_{13} + x_5(1-x_7),\\
  &-x_{12} - x_4(1-x_7), x_1^2-x_7(1-x_7),  x_6(1-x_7)),\\
 (JE_2)_x =\,& (-x_{26}-x_3(1-x_7), x_{25},  x_{24} + x_1(1-x_7), -x_{23} + x_6(1-x_7), \\
  &-x_{2}^2 + x_7(1-x_7), x_{12} - x_4(1-x_7), -x_5(1-x_7)),\\
 (JE_3)_x =\,& (-x_{36}+x_2(1-x_7), x_{35}-x_1(1-x_7), x_{34}, -x_3^2 + x_7(1-x_7),\\
  &-x_{23} - x_6(1-x_7), x_{13} + x_5(1-x_7),  -x_4(1-x_7)),\\
 (JE_4)_x =\,& (-x_{46}-x_5(1-x_7), x_{45}-x_6(1-x_7), x_4^2 - x_7(1-x_7), -x_{34},\\
  &-x_{24} + x_1(1-x_7), x_{14} + x_2(1-x_7),  x_3(1-x_7)),\\
 (JE_5)_x =\,& (-x_{56}+x_4(1-x_7), x_{5}^2-x_7(1-x_7), x_{45} + x_6(1-x_7), \\
  &-x_{35}-x_1(1-x_7), -x_{25}, x_{15} - x_3(1-x_7),  x_2(1-x_7)),\\
 (JE_6)_x =\,& (-x_{6}^2+x_7(1-x_7), x_{56}+x_4(1-x_7), x_{46} -x_5(1-x_7), -x_{36} - x_2(1-x_7), \\
  &-x_{26} + x_3(1-x_7), x_{16},   - x_1(1-x_7)).
  \end{align*}

Now  we take the natural metric $g$ on $S^6-\{p\}$. Thus,
$(S^6-\{p\}, g, J)$ is a nearly K\"ahler manifold and hence has an
$SU(3)$-structure. The K\"ahler form, $F(X,Y)=g(JX,Y)$, for any $X,
Y$ vector fields on $S^6-\{p\}$, has the form
  $$
  F=\left(\sum_{i=1}^6 x_idx_i\right) \wedge (-\beta) /(1-x_7) + \beta_1 + \beta
  \wedge dx_7 / (1-x_7),
  $$
where $\beta$ is the $1$-form
  $$
  \beta= x_6 dx_1 - x_1 dx_6 + x_2 dx_5 - x_5 dx_2 - x_4 dx_3 + x_3
  dx_4,
  $$
and $\beta_1$ is the $2$-form given by
  \begin{align*}
  \beta_1=\, & x_7 (-dx_{16}+dx_{25}+dx_{34}) + x_{1} dx_{23} +
  x_{3} dx_{12} - x_{2} dx_{13} +  x_{1} dx_{45}\\ & + x_{5} dx_{14}
  - x_{4} dx_{15} +  x_{2} dx_{46} + x_{6} dx_{24} - x_{4} dx_{26} -
  x_{3} dx_{56} - x_{6} dx_{35} + x_{5} dx_{36}.
  \end{align*}
Now, using that $ \sum_{i=1}^7 x_idx_i=0$, it follows that
  $$
  F= \beta\wedge dx_7 + \beta_1.
  $$
Then it is easy to obtain
  $$
  dF= 3 (dx_{257} + dx_{347} -dx_{167} + dx_{123} +  dx_{145}  +
  dx_{246}- dx_{356}).
  $$
Now, a long calculation shows that  $JdF$ is expressed, in terms of
the $\alpha_{ijk}$, as
  \begin{align*}
  JdF =\, & 3(1-x_7)^2 (x_{1}^2 +x_{2}^2 +x_{3}^2 +x_{4}^2 +x_{5}^2 +x_{6}^2 +x_{7}^2 )
  \Big( (-x_{34} - x_{25} + x_{16}) \alpha_{123}\\
  &- (x_{1}^2 + x_{2}^2 + x_{4}^2  - x_{7} (1- x_{7}))\alpha_{124}
  +(x_{23} - x_{45} - x_6 (1-x_7)) \alpha_{125} \\
  &- (x_{13} + x_{46} - x_5 (1-x_7))\alpha_{126}
  - (x_{23} - x_{45} + x_6 (1-x_7)) \alpha_{134} \\
  &+ (x_{1}^2 + x_{3}^2 + x_{5}^2  - x_{7} (1- x_{7}))\alpha_{135}
  +(x_{12} + x_{56} + x_4 (1-x_7)) \alpha_{136} \\
  &+ (x_{34} + x_{25} + x_{16} )\alpha_{145} - (x_{15} - x_{26} + x_3 (1-x_7)) \alpha_{146} \\
  &+ (x_{14} - x_{36} - x_2 (1-x_7))\alpha_{156}
  +(x_{13} + x_{46} + x_5 (1-x_7)) \alpha_{234} \\
  &+ (x_{12} + x_{56} - x_4 (1-x_7)) \alpha_{235}
  +(x_{2}^2 + x_{3}^2 + x_{6}^2  - x_{7} (1- x_{7}))\alpha_{236} \\
  &- (x_{15} - x_{26} - x_3 (1-x_7)) \alpha_{245}
  - (-x_{34} + x_{25} + x_{16}) \alpha_{246} \\
  &- (-x_{24} - x_{35} - x_1 (1-x_7)) \alpha_{256}
  - (x_{14} - x_{36} + x_2 (1-x_7)) \alpha_{345} \\
  &- (x_{24} + x_{35} - x_1 (1-x_7)) \alpha_{346}
  +(x_{34} - x_{25} + x_{16}) \alpha_{356} \\
  &+ (x_{4}^2 + x_{5}^2 + x_{6}^2  - x_{7} (1- x_{7}))\alpha_{456}\Big).
  \end{align*}
The $3$-forms $\Psi_{+}$ and  $\Psi_{-}$ of the $SU(3)$-structure on
$S^6-\{p\}$  are given by
  \begin{align*}
  \Psi_{+} =\, & \frac{1}{3} dF=dx_{257} + dx_{347} -dx_{167} +
  dx_{123} + dx_{145}  +
  dx_{246}- dx_{356},\\
  \Psi_{-} =\, & \frac{1}{3} JdF = \frac{1}{(1-x_7)}\big( -x_4
  dx_{127} + x_5 dx_{137}\\ &+ x_2 dx_{147} - x_3 dx_{157} + x_6
  dx_{237} - x_1 dx_{247}
  -x_3 dx_{267} + x_1 dx_{357} \\
  &+ x_2 dx_{367} + x_6 dx_{457} - x_5 dx_{467} + x_4 dx_{567}\big)
  + {\rm terms \hspace{1mm} not \hspace{1mm} containing}\hspace{1mm}
  dx_{7}.
  \end{align*}

\subsubsection{The $SU(2)$-structure on $S^5$}

Let us consider $S^5 = \{(x_1,\cdots,x_6) \in \RR^6 \mid
\sum_{i=1}^6 x_{i}^2 =1\} \subset S^6$, and $\mathbb N =
\frac{\partial}{\partial x_7}$ the unit normal vector field to
$S^5$. Then, using \eqref{hyp1}, the $SU(2)$-structure $(\eta,
\omega_i)$ on $S^5$ is given by
  \begin{equation}\label{s5hypo}
  \begin{aligned}
  \eta= & - \frac{\partial}{\partial x_7} \lrcorner F = x_6 dx_1- x_1dx_6 +
  x_2dx_5 -x_5 dx_2 +x_3 dx_4 - x_4dx_3,\\
  \omega_1 = & \, f^*(F) = x_3dx_{12}-x_2dx_{13}+x_1dx_{23}+ x_5dx_{14}-x_4dx_{15}+x_1dx_{45} \\
  &\qquad \qquad +x_6 dx_{24}-x_4dx_{26}+x_2dx_{46} +x_5dx_{36}-x_6dx_{35} -x_3dx_{56}, \\
  \omega_2 = & \, \frac{\partial}{\partial x_7} \lrcorner \Psi_{-}
  = -x_4 dx_{12} + x_5 dx_{13} + x_2 dx_{14} - x_3 dx_{15} + x_6 dx_{23} - x_1
  dx_{24}\\
  & \qquad \qquad \quad -x_3 dx_{26} + x_1 dx_{35} + x_2 dx_{36} + x_6 dx_{45} - x_5 dx_{46}
  + x_4 dx_{56},\\
  \omega_3 = &- \frac{\partial}{\partial x_7} \lrcorner \Psi_{+}   = dx_{16} - dx_{34} -  dx_{25}.
  \end{aligned}
  \end{equation}

Next we show that the $SU(2)$-structure on  $S^5$ defined by
\eqref{s5hypo} satisfies Lemma \ref{lem1}. First, we see that
  $$
  d\eta = -2(dx_{16} - dx_{34} -  dx_{25}) = -2 \omega_3.
  $$
The expression of $\omega_1$ gives
  $$
  d\omega_1=  3(dx_{123} +dx_{246} +dx_{145}-dx_{356}).
  $$
Using that $\sum_{i=1}^6 x_{i}^2=1$, so $\sum_{i=1}^6 x_{i}dx_i=0$
on $S^5$, we verify
  \begin{align*}
  3\eta\wedge\omega_2 = d\omega_1,\qquad
  -3\eta\wedge\omega_1 =  d\omega_2
  \end{align*}

Now, we apply Theorem~\ref{newnk} and Theorem~\ref{newg2} to get

\begin{thm}\label{s5nk}
Let $(S^5,\eta,\omega_i, g_{5})$ be the standard Sasaki-Einstein
5-sphere endowed with the $SU(2)$-structure determined by
\eqref{s5hypo}. Then
\begin{itemize}
  \item[i)] The $SU(3)$-structure on $S^5\times [0,\pi]$ defined by
  \eqref{newnks} is a nearly K\"ahler structure generating the round metric on the 6-sphere,
  $g_6=dt^2+\sin^2t\ g_{5}$, with two conical singularities at $t=0,
  t=\pi$.
  \item[ii)] The $G_2$-structure on $(S^5\times [0,\pi])\times [0,\pi]$
  defined by \eqref{oldg2} is a nearly parallel $G_2$-structure
  generating the roumd metric on the 7-sphere, $g_7=dq^2+\sin^2q(dt^2+\sin^2t\ g_{5})$ with
  singularities at $t=0, t=\pi, q=0, q=\pi$.
\end{itemize}
\end{thm}

\subsection{The Nearly K\"ahler structure on $S^2\times S^3\times\RR$}
As in the previous example, first we describe explicitly

\subsubsection{The standard $SU(3)$-structure on $S^3\times S^3$}
Let us consider the sphere $S^3$, viewed as the Lie group $SU(2)$,
with the basis of left-invariant $1$-forms
$\{\alpha_1,\alpha_2,\alpha_3\}$ satisfying
  \begin{gather*}
  d\alpha_1=-\alpha_2\wedge\alpha_3, \quad
  d\alpha_2=\alpha_1\wedge\alpha_3, \quad
  d\alpha_3=-\alpha_1\wedge\alpha_2.
  \end{gather*}
Denote by $\{\beta_1,\beta_2,\beta_3\}$ another basis on a second
sphere $S^3$ satisfying the same relations. Then, a nearly K\"ahler
structure on $S^3\times S^3$ is given by (\cite{ADHL}, \cite{BS})
  \begin{gather*}
  F= \frac{i}{2}(\mu_1\wedge\overline{\mu}_1 +
  \mu_2\wedge\overline{\mu}_2 +
  \mu_3\wedge\overline{\mu}_3),\quad\quad
  \Psi=i(\mu_1\wedge\mu_2\wedge\mu_3),
  \end{gather*}
where $\mu_j= \frac13 (\alpha_j + e^{\frac{2\pi i}{3}} \beta_j)$,
for $j=1,2,3$.

In terms of the real forms $\{\alpha_j,\beta_j\}$, the forms $F$,
$\Psi_+$ and $\Psi_-$ are expressed as
  \begin{align*}
  F=\, &{\sqrt{3} \over 18} (\alpha_1\wedge\beta_1 +
  \alpha_2\wedge\beta_2+\alpha_3\wedge\beta_3),\\
  \Psi_+=\, &{\sqrt{3} \over 54} (-\alpha_{12}\wedge\beta_3 +
  \alpha_{13}\wedge\beta_2-\alpha_{23}\wedge\beta_1 +
  \alpha_{1}\wedge\beta_{23} -
  \alpha_{2}\wedge\beta_{13}+ \alpha_{3}\wedge\beta_{12}),\\
  \Psi_-=\, &{1 \over 54} (2\alpha_{123}-\alpha_{12}\wedge\beta_3 +
  \alpha_{13}\wedge\beta_2-\alpha_{23}\wedge\beta_1 -
  \alpha_{1}\wedge\beta_{23} + \alpha_{2}\wedge\beta_{13}-
  \alpha_{3}\wedge\beta_{12} + 2\beta_{123}).
  \end{align*}
It is easy to check that the corresponding metric on $S^3\times S^3$
is
  \begin{gather}
  g=\frac19(\alpha_1^2+\alpha_2^2+\alpha_3^2+
  \beta_1^2+\beta_2^2+\beta_3^2 -\alpha_1\beta_1 - \alpha_2\beta_2
  -\alpha_3\beta_3).\label{metric}
  \end{gather}

\subsubsection{The $SU(2)$-structure on $S^2\times S^3$}
In order to show explicitly the induced $SU(2)$-structure on the
hypersurface $S^2\times S^3$, we first describe $S^3\times S^3$ as
the submanifold of $\RR^8$,
  \begin{gather*}
  S^3\times S^3 = \{ (x_1,\ldots,x_4,x_5,\ldots,x_8)\in \RR^8 \mid
  x_1^2+\cdots +x_4^2=x_5^2+\cdots +x_8^2=1 \}.
  \end{gather*}
With this description, we can identify
  \begin{align*}
  &\alpha_1= 2x_4dx_1 +2x_3dx_2 - 2x_2dx_3 - 2x_1dx_4, \quad\quad\
  \beta_1= 2x_8dx_5 + 2x_7dx_6- 2x_6dx_7- 2x_5dx_8,\\
  &\alpha_2= -2x_3dx_1 +2x_4dx_2+ 2x_1dx_3- 2x_2dx_4,\quad\ \,
  \beta_2= -2x_7dx_5 +2x_8dx_6+ 2x_5dx_7- 2x_6dx_8,\\
  &\alpha_3=2x_2dx_1 -2x_1dx_2+ 2x_4dx_3- 2x_3dx_4,\quad\quad\
  \beta_3=2x_6dx_5 -2x_5dx_6+ 2x_8dx_7- 2x_7dx_8.
  \end{align*}
We shall denote by $\{U_j,V_j\}_{j=1}^3$ the basis of vector fields
on $S^3\times S^3$ dual to $\{\alpha_j, \beta_j\}_{j=1}^3$.

Let us consider the hypersurface $S^2\times S^3 \subset S^3\times
S^3$ given by $x_4=0$. Then, with respect to the metric
\eqref{metric}, the vector field
  \begin{gather*}
  \mathbb N = -\sqrt{3} (2 x_1 U_1 + 2x_2 U_2 +2x_3 U_3 + x_1 V_1 +
  x_2 V_2 + x_3 V_3)
  \end{gather*}
is a unit normal vector field along $S^2\times S^3$.

Next, we describe explicitly the induced $SU(2)$-structure
\eqref{hyp1}, taking $f$ as the inclusion map.

A direct calculation, using that $x_1 \alpha_1 + x_2 \alpha_2 +
x_3\alpha_3\equiv 0$ on $S^2\times S^3$, shows that the form $\eta$
is expressed as
  \begin{equation}\label{e-s-3}
  \begin{aligned}
  \eta=\, &-\mathbb N\lrcorner F = {1\over 3} (x_1\beta_1 +x_2\beta_2 +
  x_3\beta_3) ={2\over 3} \big((x_{18}-x_{27}+x_{36})dx_5 \\
  &+  (x_{17}+x_{28}-x_{35})dx_6 + (-x_{16}+x_{25}+x_{38})dx_7
  +(-x_{15}-x_{26}-x_{37})dx_8 \big).
  \end{aligned}
  \end{equation}

Since $f$ is the inclusion, taking $x_4=0$ in the expressions of
$\alpha_j$ above, we get
  \begin{align*}
  \omega_1=f^*F =\, & \frac{2\sqrt{3}}{9} \Big((x_{26}+x_{37})dx_{15} +
  (-x_{25}-x_{38})dx_{16} + (x_{28}-x_{35})dx_{17} \\
  &+(-x_{27}+x_{36})dx_{18} +
  (-x_{16}+x_{38})dx_{25} + (x_{15}+x_{37})dx_{26} \\
  &+(-x_{18}-x_{36})dx_{27} +
  (x_{17}-x_{35})dx_{28} + (-x_{17}-x_{28})dx_{35} \\
  &+(x_{18}-x_{27})dx_{36} + (x_{15}+x_{26})dx_{37} +
  (-x_{16}+x_{25})dx_{38}\Big).
  \end{align*}

For computing $\omega_2$ and $\omega_3$, take into account the
equality $x_3\alpha_1\wedge\alpha_2 - x_2\alpha_1\wedge\alpha_3 +
x_1\alpha_2\wedge\alpha_3 = 4 (x_3dx_{12} -x_2 dx_{13}+
x_1dx_{23})$, to get
  \begin{align*}
  \omega_2=\mathbb N\lrcorner\Psi_-
  =\, & \frac{2\sqrt{3}}{9} \Big( -x_3dx_{12} +x_2 dx_{13} - x_1dx_{23} + x_8 dx_{15}
  +x_7dx_{16} -x_6dx_{17} - x_5dx_{18} \\
  &-x_7dx_{25} + x_8dx_{26} + x_5dx_{27} -x_6dx_{28} + x_6dx_{35} -
  x_5dx_{36} + x_8dx_{37} - x_7dx_{38}\Big);\\
  \omega_3=-\mathbb N\lrcorner\Psi_+
  =\, & \frac{2}{9} \Big(-x_3dx_{12} +x_2 dx_{13} - x_1dx_{23} - x_8 dx_{15}
  -x_7dx_{16} +x_6dx_{17} + x_5dx_{18} \\
  &+x_7dx_{25} - x_8dx_{26} - x_5dx_{27} +x_6dx_{28} - x_6dx_{35} +
  x_5dx_{36} - x_8dx_{37} + x_7dx_{38}\Big)\\
  &+ \frac{4}{9} \Big(x_3dx_{56} -x_2 dx_{57} + x_1 dx_{58} +x_1dx_{67}
  +x_2dx_{68} + x_3dx_{78}\Big).
  \end{align*}

Notice that $S^2\times S^3$ is {\em not} a totally geodesic
hypersurface of $S^3\times S^3$; for example, for $T=x_2U_1-x_1U_2$
which is tangent to $S^2\times S^3$, we have
  $$
  g(\nabla_{T} \mathbb N,V_3) = -\frac{\sqrt{3}}{36} (x_1^2+x_2^2),
  $$
which is non-zero on $S^2\times S^3$, and thus the second
fundamental form does not vanish identically. Therefore, we cannot
apply Lemma \ref{lem1} to establish that the $SU(2)$-structure
$(\eta,\omega_i)$ induced on $S^2\times S^3$ from the nearly
K\"ahler structure of $S^3\times S^3$ is hypo. To solve this
problem, we proceed as follows. We have
  \begin{align*}
  d\eta= \, & \frac{1}{3}(dx_{1}\wedge\beta_1-x_{1}\wedge\beta_{23}+dx_{2}\wedge\beta_{2}+
x_{2}\wedge\beta_{13}+dx_{3}\wedge\beta_{3}-x_{3}\wedge\beta_{12})\\
= \, &\frac{2}{3}\Big(x_8 dx_{15}+x_7dx_{16} -x_6dx_{17} - x_5dx_{18} -x_7dx_{25}\\
  &+ x_8dx_{26} + x_5dx_{27} -x_6dx_{28} + x_6dx_{35} -
  x_5dx_{36} + x_8dx_{37} - x_7dx_{38}\\
  &-2x_3dx_{56} +2x_2 dx_{57} -2 x_1 dx_{58} -2x_1dx_{67} -2x_2dx_{68}
  -2 x_3dx_{78} \Big),
  \end{align*}
and so we  can write
  $$
  \omega_3=-\frac{1}{3}d\eta + \frac{2}{9}(-x_3dx_{12} +
  x_2dx_{13}-x_1dx_{23}),
  $$
which implies that
  $$
  d\eta\not=-2\omega_3,
  $$
since the form $-x_3dx_{12} + x_2dx_{13}-x_1dx_{23}$ is the standard
volume form on $S^2$, and
  $$
  d\omega_3=0,
  $$
because $d(-x_3dx_{12} + x_2dx_{13}-x_1dx_{23})=0$ on $S^2\times
S^3$. Moreover, we get
  \begin{gather*}
  d\omega_1 = 3 \eta\wedge \omega_2, \\
  d\omega_2 = -3 \eta\wedge \omega_1.
  \end{gather*}

The previous equalities show that $(\eta, \omega_i)$ is hypo on
$S^2\times S^3$, but it does not satisfy equations \eqref{e-s}
because $d\eta\not=-2\omega_3$.

On the other hand, a direct calculation shows that
  $$
  \eta\wedge (d\eta)^2= -\frac{2}{27}(x_3dx_{12} -x_2 dx_{13}+
  x_1dx_{23})\wedge\beta_{123} \not= 0,
  $$
so $\eta$ is a contact form on $S^2\times S^3$.

\begin{rem}
Let us see that $3\eta=x_1\beta_1+x_2\beta_2+x_3\beta_3 \in
\Omega^1(S^2\times S^3)$ is the natural contact form on $S^2\times
S^3$ seen as the tangent sphere bundle over $S^3$ (see
\cite{Blair}). As $S^3$ is parallelizable, the tangent bundle to
$S^3$ is isomorphic to $\RR^3\times S^3$. Let $V_1,V_2,V_3$ be an
orthonormal basis of left-invariant vector fields, and let
$\beta_1,\beta_2,\beta_3$ be the dual basis of left-invariant
$1$-forms. The isomorphism $\RR^3\times S^3 \cong TS^3$ is given by
$((a_1,a_2,a_3),p) \mapsto \sum a_i V_i(p)$. The metric of $S^3$ is
$g=\beta_1^2 +\beta_2^2 +\beta_3^2$. Consider the unit sphere in the
tangent bundle $T_1S^3 \cong S^2\times S^3$. If $x_1,x_2,x_3$ are
the natural coordinates in the $\RR^3$ factor of $\RR^3\times S^3$,
then $T_1S^3$ is given by the equation $x_1^2+x_2^2+x_3^2=1$.

The natural $1$-form of $T^*S^3$ (the Liouville form) is given as
$\lambda\in \Omega^1(T^*S^3)$, $\lambda_\alpha(v)= \alpha(d\pi(v))$,
where $\pi:T^*S^3\to S^3$. Using the metric, we identify
$g:TS^3\cong T^*S^3$. Then $g^*\lambda|_{T_1S^3}$ is the natural
contact form for $T_1S^3$ (see \cite{Blair}).

It is easy to see that $3\eta=g^*\lambda|_{T_1S^3}$. Actually,
$x_1\beta_1+x_2\beta_2+x_3\beta_3=g^*\lambda \in \Omega^1(TS^3)$.
Equivalently, we need to see that
$y_1\beta_1+y_2\beta_2+y_3\beta_3=\lambda \in \Omega^1(T^*S^3)$,
where $y_1,y_2,y_3$ are the coordinates of the $\RR^3$ factor of
$T^*S^3\cong \RR^3\times S^3$. But take $\alpha= \sum
a_i\beta_i(p)\in T^*_pS^3$. Then
 $$
 \lambda_\alpha(v_1,v_2)=\alpha(v_2)= \sum a_i \beta_i(p)(v_2)=
  \left(\sum y_i \beta_i\right) ((a_1,a_2,a_3),p)(v_1,v_2),
 $$
for $(v_1,v_2)\in T_\alpha(\RR^3\times S^3)=T_\alpha(T^*S^3)$,
identifying $\beta_i$ in $S^3$ with its pull-back to $\RR^3\times
S^3$.
\end{rem}

The following result shows how the hypo structure on $S^2\times S^3$
described above can be deformed into a double hypo structure, and
even into a Sasaki-Einstein structure.

\begin{pro}\label{familyS2xS3}
Let $(\eta, \omega_i)$ be the hypo structure on $S^2\times S^3$
given above. For each $\lambda<0$ and $\mu > \frac{\lambda}{3}$, the
quadruplet
\begin{equation}\label{e-s-1-bis}
(\, \tilde\eta=\eta,\quad
\tilde\omega_1=\sqrt{3\lambda(\lambda-3\mu)}\, \omega_1,\quad
\tilde\omega_2=\sqrt{3\lambda(\lambda-3\mu)}\, \omega_2,\quad
\tilde\omega_3= \lambda\, d\eta + \mu\, vol_{S^2}\,)
\end{equation}
defines a hypo structure on $S^2\times S^3$, which is double hypo if
and only if $\lambda<-\frac{1}{4}$ and $\mu=\frac{\lambda(4\lambda
+2)}{3(4\lambda +1)}$. Moreover, the
$SU(2)$-structure~(\ref{e-s-1-bis}) is Sasaki-Einstein only for
$\lambda=-\frac{1}{2}$ and $\mu=0$.
\end{pro}

\begin{proof}
Since $(\eta,\omega_i)$ is a SU(2)-structure and $d\eta\wedge
d\eta=-\frac{2}{3} d\eta\wedge vol_{S^2}$, we have that
$$\tilde\omega_i\wedge \tilde\omega_i= \lambda(\lambda-3\mu)
d\eta\wedge d\eta
$$
for $i=1,2,3$. Moreover, $\omega_i\wedge vol_{S^2}=0$ and
$\omega_i\wedge d\eta=0$ for $i=1,2$, so the
quadruplet~\eqref{e-s-1-bis} satisfies~\eqref{defsu2}.

In order to see that~\eqref{e-s-1-bis} also satisfies
condition~\eqref{defsu2-2}, let $X=\sum_{i=1}^3 (f_i U_i + a_i V_i),
Y=\sum_{i=1}^3 (g_i U_i + b_i V_i)$ be vector fields such that
$X\lrcorner\tilde\omega_1=Y\lrcorner\tilde\omega_2$. This condition
implies that
$$
\begin{array}{rl}
x_2 g_3 -x_3 g_2 =\!\!\!&\!\! x_3(x_3f_1-x_1f_3)-x_2(x_1f_2-x_2f_1),\\
x_3 g_1 -x_1 g_3= \!\!\!&\!\! x_1(x_1f_2-x_2f_1)-x_3(x_2f_3-x_3f_2),\\
x_1 g_2 -x_2 g_1= \!\!\!&\!\! x_2(x_2f_3-x_3f_2)-x_1(x_3f_1-x_1f_3),
\end{array}
$$
and
$$
\begin{array}{rl}
b_1=\!\!\!&\!\! x_2 a_3 -x_3 a_2 + x_3(x_3 g_1 -x_1 g_3)-x_2(x_1 g_2
-x_2 g_1),\\
b_2= \!\!\!&\!\! x_3 a_1 -x_1 a_3 + x_1(x_1 g_2 -x_2 g_1)-x_3(x_2
g_3 -x_3 g_2),\\
b_3= \!\!\!&\!\! x_1 a_2 -x_2 a_1 + x_2(x_2 g_3 -x_3 g_2)-x_1(x_3
g_1 -x_1 g_3).
\end{array}
$$
Then, on $S^2\times S^3$ we have that
$$
d\eta(X,Y)=-\frac{1}{6} \left(b_1^2+b_2^2+b_3^2+(x_2 a_3 -x_3
a_2)^2+(x_3 a_1 -x_1 a_3)^2+(x_1 a_2 -x_2 a_1)^2 \right)
$$
and
$$
vol_{S^2}(X,Y)= \frac{1}{4} \left( (b_1-x_2 a_3 +x_3 a_2)^2+(b_2-
x_3 a_1 +x_1 a_3)^2+ (b_3-x_1 a_2 +x_2 a_1 )^2\right).
$$
Therefore, $(\lambda\, d\eta+ \mu\, vol_{S^2})(X,Y) \geq 0$ when
$\lambda<0$ and $\mu>\lambda/3$.

The SU(2)-structure~(\ref{e-s-1-bis}) clearly satisfies that
$d\tilde\omega_1 = 3 \tilde\eta\wedge \tilde\omega_2$,
$d\tilde\omega_2 = -3 \tilde\eta\wedge \tilde\omega_3$ and
$d\tilde\omega_3=0$. Moreover, using again that $d\eta\wedge
d\eta=-\frac{2}{3} d\eta\wedge vol_{S^2}$, the structure is double
hypo if and only if $2 \lambda^2+ \lambda -6\lambda\, \mu-
\frac{3}{2}\mu =0$. Therefore, $\lambda\not=-\frac{1}{4}$ and $\mu=
\frac{\lambda}{3}(1+\frac{1}{4\lambda+1})$ in order the latter
relation be satisfied. Since $\mu>\lambda/3$, we must have
$\lambda<-1/4$. Finally, the $SU(2)$-structure~\eqref{e-s-1-bis}
satisfies equations \eqref{e-s}, i.e. it is a Sasaki-Einstein hypo
structure, only for $\lambda=-1/2$ and $\mu=0$.
\end{proof}

Now, we apply Proposition~\ref{familyS2xS3} and Theorems~\ref{newnk}
and~\ref{newg2} to get
\begin{thm}\label{s6nk}
Let $(S^2\times S^3,\tilde\eta,\tilde\omega_i, g)$ be the
Sasaki-Einstein manifold endowed with the $SU(2)$-structure
determined by \eqref{e-s-1-bis} for $\lambda=-1/2$ and $\mu=0$. Then
\begin{itemize}
  \item[i)] The $SU(3)$-structure on $S^2\times S^3\times [0,\pi]$ defined by
  \eqref{newnks} is a nearly K\"ahler structure generating the metric
  $g_6=dt^2+\sin^2t\ g$ with two conical singularities at $t=0, t=\pi$.
  \item[ii)] The $G_2$-structure on $(S^2\times S^3\times S^1)\times
  [0,\pi]$ defined by \eqref{oldg2} is a nearly parallel $G_2$-structure
  generating the metric $g_7=dq^2+\sin^2q(dt^2+\sin^2t\ g)$ with
  singularities at $t=0, t=\pi, q=0, q=\pi$.
\end{itemize}
\end{thm}

\subsection{The Nearly K\"ahler structure on $Y^{p,q}\times \RR$}
We start with the recently discovered in \cite{GMS} infinite family
of Sasaki-Einstein metric on $S^2\times S^3$, labeled by two coprime
integers $p>1, q<p$ and refered as $Y^{p,q}$. Geometrically they are
all $U(1)$-bundles over an axially squashed $S^2$ bundle over a round
$S^2$. We take the explicit local description of the Sasaki-Einstein
$SU(2)$-structure presented in \cite{MS}.  In terms of local
coordinates $y,\beta,\theta,\phi,\psi$ they can be described as
follows:
\begin{gather}
\eta=\frac13\left(d\psi-\cos\theta d\phi+y(d\beta+c\cos\theta d\phi)\right),\nonumber\\
\omega_3=\frac16\left((cy-1)\sin\theta d\theta\wedge
d\phi-dy\wedge(d\beta+c\cos\theta
d\phi)\right),\nonumber\\\label{newpq}
\omega_2=\sqrt{\frac{1-cy}{6w(y)r(y)}}\left(d\theta\wedge
dy-\frac{w(y)r(y)\sin\theta}{6}d\theta\wedge
d\beta\right),\\\nonumber
\omega_1=\sqrt{\frac{1-cy}{6w(y)r(y)}}\left(\sin\theta d\phi\wedge
dy+\frac{w(y)r(y)}{6}d\theta\wedge(d\beta+c\cos\theta d\phi)\right),
\end{gather}
where $$w(y)=\frac{2(a-y^2)}{1-cy}, \qquad
r(y)=\frac{a-3y^2+2cy^3}{a-y^2},$$ and $a,c$ are constants. If $c=0$
one can obtain the known homogeneous metric on $S^2\times S^3$ and
for $c=1=a$ one can recover the round metric on 5-sphere $S^5$.
However, for $c\not=0, 0<a<1$ one can get  irregular Sasaki-Einstein
structures, i.e. the orbits of the Killing vector field dual to
$\eta$ are non-compact \cite{GMS}.

It is easy to check that the $SU(2)$-structure \eqref{newpq} satisfies
\eqref{e-s}. Apply Theorem~\ref{newnk} and Theorem~\ref{newg2} to
get
\begin{thm}\label{pqnk}
Let $(Y^{p,q},\eta,\omega_i, g)$ be the Sasaki-Einstein manifold
endowed with the $SU(2)$-structure determined by \eqref{newpq}. Then
\begin{itemize}
  \item[i)] The $SU(3)$-structure on $Y^{p.q}\times [0,\pi]$ defined by
  \eqref{newnks} is a nearly K\"ahler structure generating the metric
  $g_6=dt^2+\sin^2t\ g$ with two conical singularities at $t=0, t=\pi$.
  \item[ii)] The $G_2$-structure on $(Y^{p,q}\times [0,\pi])\times
  [0,\pi]$ defined by \eqref{oldg2} is a nearly parallel $G_2$-structure
  generating the metric $g_7=dq^2+\sin^2q(dt^2+\sin^2t\ g)$ with
  singularities at $t=0, t=\pi, q=0, q=\pi$.
\end{itemize}
\end{thm}
\subsection{Evolution which is not an $SU(2)$-structure}\label{last}

For half flat $SU(3)$-structures Hitchin shows \cite{Hit} that if his
evolution equations are satisfied, and for $t=0$ the structure is
half-flat, then the half flat $SU(3)$ condition is preserved in time
provided some non-degeneracy condition for the evolved
$SU(3)$-structure holds.

For a hypo and nearly hypo $SU(2)$-structure we find an example which
solves the Conti-Salamon and our nearly hypo evolution equations but
there exists a solution to the evolution equations which is not an
$SU(2)$-structure, i.e. the situation is a little bit different.

We take the double hypo structure on the Lie group isomorphic to
$SU(2)\times A^2$ defined in Proposition~\ref{red-equations} by
\eqref{structure} and \eqref{differentials} for $\mu=0$.

We find the following solution to the Conti-Salamon hypo evolution
equations \eqref{vol}
$$\eta(t)=\eta,\quad \omega_1(t)=\omega_1-t\,d\eta,\quad
\omega_3(t)=-\sinh3t\,\omega_1+\omega_3,\quad
\omega_2(t)=\cosh3t\,\omega_2$$ which is not an SU(2)-structure
since $\omega_2^2=\omega_3^2=\cosh^23t$ while $\omega_1^2=1$.

We obtain the following solution to the nearly-hypo evolution
equations \eqref{evolunk}
\begin{gather*}
\eta(t)=\eta, \quad \omega_2(t)=\cos\sqrt 3 \,t\,\omega_2,\\
\omega_1(t)=\cos\sqrt 3 \,t\,\omega_1 -\frac{\sqrt 3}{2}\sin2\sqrt 3
\,t\,e^{14}+\frac{1}{2\sqrt 3}\sin2\sqrt 3 \,t\,e^{23},\\
\omega_3(t)=\frac1{\sqrt 3}\sin\sqrt 3 \,t\,\omega_1+\cos2\sqrt 3
\,t\,e^{14}+\left(\frac43-\frac13\cos2\sqrt 3 \,t\right)e^{23}
\end{gather*}
which is not an $SU(2)$-structure again.

\medskip
\noindent {\bf Acknowledgments.}  This work has been partially
supported through grants MCyT (Spain) MTM2004-07090-C03-01,
MTM2005-08757-C04-02 and Project UPV 00127.310-E-15909/2004. We
would like to thank Simon Salamon for stimulating discussions on the
existence of nearly parallel $G_2$-structures with conical
singularities, and Diego Conti for informing us on his answer to
Question $1$. We also thank Anna Fino for helpful discussions, and
the referee for useful comments. In particular, the referee
suggested the inclusion of diagrams
which improved greatly the exposition.

\bibliographystyle{hamsplain}

\begin{thebibliography}{19}
\bibitem{ADHL}
B.S. Acharya, F. Denef, C. Hofman, N. Lambert, {\em Freund-Rubin
Revisited}. Preprint {\tt hep-th/0308046}.

\bibitem{BFGK}
H. Baum, Th. Friedrich, R. Grunewald, I. Kath, {\em Twistors and
Killing spinors on Riemannian manifolds}, Seminarbericht Nr. 108,
Humboldt-Universit\"at zu berlin, 1990.

\bibitem{BM}
A. Bilal, S. Metzger, {\em Compact weak $G_2$-manifolds with conical
singularities},  Nuclear Phys. {\bf B 663} (2003), no. 1-2,
343--364.

\bibitem{Blair}
D.E. Blair, Riemannian geometry of contact and symplectic manifolds,
{\em Progress in Mathematics 203}, Birkh\"{a}user Boston, Inc., Boston,
MA, 2002.

\bibitem{BGal}
C.Boyer, K. Galicki, {\em 3-Sasakian manifolds}, Surveys Diff. Geom.
{\bf 6} (1999), 123--184.

\bibitem{BG}
C.Boyer, K. Galicki, {\em New Einstein metrics in dimension five},
J. Diff. Geom. {\bf 57} (2001), 443--463.

\bibitem{BGN}
C.Boyer, K. Galicki, M. Nakamae,
 {\em On the geometry of Sasakian-Einstein 5-manifolds},
 Math. Ann.  {\bf 325} (2003), 485--524.

\bibitem{Br1}
R.L. Bryant, {\em Submanifolds and special structures on the
octonians}, J. Differential Geom. {\bf 17} (1982), 185--232.

\bibitem{BS} R. Bryant, S.Salamon, {\em On the construction of some
  complete metrics with exceptional holonomy}, Duke Math. J. {\bf 58}
  (1989), 829-850.

\bibitem{Cab1}
F.M. Cabrera, {\em $SU(3)$-structures on hypersurfaces of manifolds
with $G_2$-structure}, Monatsh. Math. {\bf 148} (2006), 29-50.

\bibitem{Cal}
E. Calabi, {\em Construction and properties of some $6$-dimensional
almost complex manifolds}, Trans. Amer. Math. Soc. {\bf 87} (1958),
407--438.

\bibitem{ConS}
D. Conti, S. Salamon, {\em Generalized Killing spinors in dimension
$5$}, Trans. Amer. Math. Soc. {\bf 359} (2007), 5319--5343.


\bibitem{Conti} D. Conti, private communication, $2008$.

\bibitem{FG}
M. Fern\'andez, A. Gray, {\em Riemannian manifolds with structure
group $G_2$}, Ann. Mat. Pura Appl.  {\bf 32}  (1982), 19--45.

\bibitem{FI1}
Th. Friedrich, S. Ivanov, {\em Parallel spinors and connections with
skew-symmetric torsion in string theory}, Asian J. Math.  {\bf 6}
(2002),  no. 2, 303--335.

\bibitem{GMS}
J. Gauntlett, D. Martelli, J. Sparks, D. Waldram, {\em
Sasaki-Einstein Metrics on $S^2\times  S^3$}, Adv.Theor. Math. Phys.
{\bf 8} (2004) 711-734.

\bibitem{Gr}
A. Gray, {\em Some examples of almost Hermitian manifolfs}, Illinois
J. Math. {\bf 10} (1969), 353--366.

\bibitem{Gr3}
A. Gray, {\em Nearly K\"ahler manifolds},  J. Differential Geometry
{\bf 4}  (1970), 283--309.

\bibitem{Gr2}
A. Gray, {\em Riemannian manifolds with geodesic symmetries of order
$3$},  J. Differential Geometry {\bf 7} (1972), 343--369.

\bibitem{Gr1}
A.Gray, {\em The structure of nearly K\"ahler manifolds},  Math.
Ann.  {\bf 223}  (1976), no. 3, 233--248.

\bibitem{Hit}
N. Hitchin, {\em Stable forms and special metrics}, Global
differential geometry: the mathematical legacy of Alfred Gray
(Bilbao, 2000),  70--89, Contemp. Math., {\bf 288}, Amer. Math.
Soc., Providence, RI, 2001.

\bibitem {Kir}
V. Kirichenko, {\em K-spaces of maximal rank}, (russian), Mat. Zam.
{\bf 22} (1977), 465-476.

\bibitem{MS} D. Martelli, J. Sparks, {\em Toric Geometry,
Sasaki-Einstein Manifolds and a New Infinite Class of AdS/CFT
Duals}, Commun. Math. Phys. {\bf 262} (2006) 51-89.

\bibitem{Stock} S. Stock, Lifting $SU(3)$-structures to nearly
parallel $G_{2}$-structures, Preprint {\tt arXiv:0707.2029v1}.

\end{thebibliography}
\providecommand{\bysame}{\leavevmode\hbox
to3em{\hrulefill}\thinspace}

\end{document}